\documentclass[conf]{new-aiaa}
\usepackage[utf8]{inputenc}

\usepackage{graphicx}
\usepackage{amsmath}
\usepackage[version=4]{mhchem}
\usepackage{siunitx}
\usepackage{longtable,tabularx}
\usepackage{dsfont}
\usepackage{calc}
\usepackage{color}
\newcommand{\norm}[1]{\left\lVert#1\right\rVert}
\newcommand{\abs}[1]{\left|#1\right|}
\newcommand{\menge}[1]{\left\lbrace#1\right\rbrace}
\newcommand{\limes}[2]{\lim\limits_{{#1}\rightarrow {#2}}}
\newcommand{\intD}{\;\mathrm{d}}
\newcommand{\rom}[1]{\uppercase\expandafter{\romannumeral #1\relax}}
\newcommand{\landauO}{\mathcal{O}}
\newtheorem{theorem}{Theorem}

\title{Windowing Regularization Techniques for Unsteady Aerodynamic Shape Optimization}

\author{
Steffen Schotthöfer\footnote{Student Researcher, sschotth@rhrk.uni-kl.de.},
Beckett Y. Zhou\footnote{Research Scientist, yuxiang.zhou@scicomp.uni-kl.de, Member AIAA.}, 
Tim Albring\footnote{Ph.D Candidate, tim.albring@scicomp.uni-kl.de.} and 
Nicolas R. Gauger\footnote{Professor, nicolas.gauger@scicomp.uni-kl.de, Associate Fellow AIAA.}
}
\affil{Chair for Scientific Computing, TU Kaiserslautern \\ Bldg 34, Paul-Ehrlich-Strasse, 67663 Kaiserslautern}

\begin{document}

\maketitle

\begin{abstract}
Unsteady Aerodynamic Shape Optimization presents new challenges in terms of sensitivity analysis of time-dependent objective functions.
In this work, we consider periodic unsteady flows governed by the URANS equations.
Hence, the resulting output functions acting as objective or constraint functions of the optimization are themselves periodic with unknown period length, that may depend on the design parameter of said optimization. Sensitivity Analysis on the time-average of a function with these properties turns out to be difficult. Therefore, we explore methods to regularize the time average of such a function with the so called windowing-approach. Furthermore, we embed these regularizers into the discrete adjoint solver for the URANS equations of the multi-physics and optimization software SU2. Finally, we exhibit a comparison study between the classical non regularized optimization procedure and the ones enhanced with regularizers of different smoothness and show that the latter result in a more robust optimization.
\end{abstract}

\section{Introduction}
Aerodynamic shape optimization has been a subject of active research for the last decades  \cite{Nadarajah2007OptimumSD,NielsenDiskin, BeckettTimNico, Albring_Sagebaum_Gauger,ThomasHallKennethDowell}. The goal of the design process usually consists of optimizing one or multiple aerodynamic coefficients such as drag or lift of the object in focus to improve its flow properties.
Since the advent of adjoint based methods \cite{Jameson1988AerodynamicDV, pironneau_1974}, for which the computational cost is independent of the number of design variables, it is possible to address many applicable large scale problems such as aerodynamic  and aerostructural optimizations of complete aircraft configurations.

Challenges of the optimization problem also arise from the nature of the equations used to model the fluid flow. Traditionally, the underlying problem is considered to be in a steady state, which is a reasonable assumption for many flow conditions.
However, many fluid flows in industry are naturally unsteady and turbulent. Such flows occur for example in aerodynamic and aeroacoustic design of rotorcraft and wind turbines or active flow control for high lift devices \cite{Nadarajah2007OptimumSD}.
In comparison to steady state aerodynamic shape optimization, the unsteady counterpart has only recently attracted more attention \cite{Nadarajah2007OptimumSD, ThomasHallKennethDowell, Albring,BeckettTimNico}. This is due to the fact that large amounts of solution data and computational time are required to solve the unsteady adjoint equation. Nevertheless, more and more research is being carried out in this area as, on the one hand, time-accurate numerical methods have been further developed and, on the other, the computing power of high-performance computing centres has increased.

For flow conditions in which large scale vortex structures are formed, which is the case, for example, in the wake of bluff bodies as an airfoil with high angle of attack, scale-resolving methods are required. Reason for this is that turbulence models are not designed for detached flows \cite{Spalart}. 
Scale resolving methods like Large Eddy Simulation (LES), Detached Eddy Simulation (DES) and Direct Numerical Simulation (DNS) resolve some part of the chaotic dynamics of turbulent fluid flows. Chaotic systems give rise to new challenges for sensitivity analysis, which is one aspect  in the solution process of the shape optimization problem. 
There are different approaches to cope with the fact that in chaotic systems, meaningful time averaged sensitivities can not be directly computed. On the one hand, finite difference methods produce wrong results, on the other hand, forward or adjoint mode of automatic differentiation produces correct results, but they are by orders of magnitude too big to be useful \cite{Albring}.
In the search for efficient algorithms for a meaningful sensitivity analysis, Least Squares Shadowing approaches were considered \cite{WangFDNILSS, WangLSS, WangNILSS}, for example.
Both scale resolving methods and sensitivity analysis of chaotic systems are computationally very expensive.

Unsteady Reynolds Averaged Navier Stokes (URANS) methods are therefore often used to calculate large scale unsteady turbulent flows. In URANS simulations the increased dissipation induced by turbulence models partially dampens the chaotic behavior of the system. In many cases, URANS governed flows exhibit limit cycle oscillations, which is the focus of this work. Here, a meaningful objective function for the shape optimization is a windowed temporal average of an aerodynamic coefficient over the length of a period of the limit cycle. A property of unsteady flows which exhibit limit cycle oscillations is that not only the shape of the limit cycles, but also the length of the period depends on the design parameters. Understanding sensitivity analysis for limit cycle oscillations in URANS-governed flows is an important first step to tackle the problem of sensitivity analysis for chaotic systems.

The goal of this publication is to apply methods of Krakos et al. \cite{Krakos}, who use a windowing smoothing approach originating from the field of signal processing to regularize the time averaged output function, in the context of aerodynamic shape optimization. 
Additionally, we present new insights into the way windowed time averaged sensitivities behave in more turbulent systems with high Reynolds number flows whose instantaneous sensitivities exhibit exponentially increasing amplitudes.
We embed the regularized objective function in the discrete optimization framework of the multi-physics CFD suite SU2 \cite{SU2_big}, where we use automatic (or algorithmic) differentiation (AD), replacing the labour-intensive and error prone manual differentiation of the discretized equations. Another advantage of AD-based adjoints is the fact that AD leads to accurate to machine precision results by construction as they do not incur any roundoff or truncation error. Furthermore, although turbulence models are not analytically differentiable, they are still  {algorithmically} differentiable. Therefore, the  {frozen turbulence} assumption, which is typically used in many URANS-based adjoint formulations, is eliminated. The AD based discrete adjoint framework inherits the same convergence properties as the primal flow solver and thus yields a robust method to compute adjoints. 

The remainder of this paper is structured as follows. In Section II, the challenges in sensitivity analysis in the presence of limit cycle oscillations are discussed and we review the windowing framework by Krakos et al. to overcome these challenges.
Section III describes the discrete primal and adjoint flow solver of SU2. Furthermore, the windowing regularization approach is embedded in the discrete adjoint solver and we show convergence properties of the adjusted adjoint solver.
In Section IV, numerical results are presented, where we validate the convergence properties of the windowed time average objective functions and we showcase where the traditional approach fails. In addition, we validate the consistency of the primal (tangent) AD mode and the adjoint AD mode for computing sensitivities of windowed time averaged objective functions. We showcase the advantages of windowed sensitivity computations in the case of high dimensional surface sensitivities compared to the traditional approach. Here, we display time averaged surface sensitivities of the NACA0012 airfoil calculated with different windows and pin down critical zones. Finally, we
perform a series of shape optimization procedures on the NACA0012-airfoil using different windows and compare the results with regard to our previous findings.
In Section V, we present numerical results for a pitching airfoil  to illustrate a test case, where the period length is not dependent on the design parameters.

\section{Computing Sensitivities of Limit Cycle Oscillations}
The sensitivity computation of a given objective function such as the time averaged drag or lift of an airfoil is an essential part of any optimization attempt. If one employs URANS for the flow simulation that acts as the optimization constraint, we arrive at some point at a spatial discretization of the form 
\begin{align}\label{eq_ODE}
	\partial_t u +R(u) &= 0,\qquad t\in\left[0,t_{f}\right]\\
	u(t=0) &= u_0
\end{align} 
where $R(U)$ is the residual obtained by a finite volume approach for solving the turbulent unsteady compressible Navier-Stokes equations, possibly with additional turbulence model. The vector $U$ denotes the solution of the URANS equations on each grid point of the computational domain and $U_0$ are the freestream conditions that display the initial values of the flow.

This approach, called method of lines, can be seen as a coupled system of ordinary differential equations that govern a dynamical system. We define the scalar instantaneous output  $g\in\mathbb{R}$ of the dynamical system, that might be the time-dependent drag $C_D$ or lift coefficient $C_L$ with additional dependence on a system parameter $\sigma\in\mathbb{R}^{n_d}$. In many applications, this dynamical system exhibits periodic behavior after some transient phase $t\in\left[0,t_{tr}\right]$. 
It makes sense to compute the mean value over a period $T(\sigma)$ as the scalar output of the dynamical system, 
\begin{align}\label{eq_tAvg}
 {J}(\sigma)=\frac{1}{T(\sigma)}\int_{t_{tr}}^{t_{tr}+T(\sigma)} g(t,\sigma)\intD t,
\end{align}
that later acts as the objective function of an optimal control problem. In the following we shift the time scale by $t_{tr}$ for the sake of easier notation.

Sensitivity computation of a limit cycle oscillator imposes several challenges, namely that one often does not know the exact duration of a period $T$ and that the period length is dependent on the system parameter $\sigma$, which might be an initial condition, a boundary condition, a design variable or an additional parameter. 
A common approach to compute the time average without prior knowledge of $T(\sigma)$ is to average over a fixed time $M$
\begin{align}
{J}_M(\sigma,M)=\frac{1}{M}\int_{0}^{M} g(t,\sigma)\intD t.
\end{align}
The corresponding sensitivity of ${J}$ and ${J}_M$ with respect to $\sigma$ is
\begin{align}
 \frac{\intD}{\intD \sigma}{J} &=  \frac{\intD}{\intD \sigma}\left(\frac{1}{T(\sigma)}\int_{0}^{T(\sigma)} g(t,\sigma)\intD t\right) \not= \frac{1}{T(\sigma)}\int_0^{T(\sigma)} \frac{\intD}{\intD \sigma} g(t,\sigma)\intD t \label{eq_avgSens1}\\
 \frac{\intD}{\intD \sigma}{J}_M &= \frac{1}{M}\int_{0}^{M}\frac{\intD}{\intD \sigma} g(t,\sigma)\intD t.
\end{align}
Note that we cannot simply interchange differentiation and integration in Eq.~\eqref{eq_avgSens1} due to dependence of $T$ on $\sigma$.
The hope is to reduce the error in $\abs{{J}_M -{J}}$ respectively $\abs{\frac{\intD}{\intD \sigma}{J}_M - \frac{\intD}{\intD \sigma}{J}}$ as $M$ grows, so the question at hand is whether $\limes{M}{\infty}{J}_M={J}$ and (more importantly) $\limes{M}{\infty}\frac{\intD}{\intD \sigma}{J}_M = \frac{\intD}{\intD \sigma}{J}$. 
To answer this question in detail,  Krakos et al. \cite{Krakos} introduced a more general notation of a weighted time average 
\begin{align}\label{eq_wndAvg1}
  {J}_w(\sigma,M) =\frac{1}{M}\int_0^M w\left(\frac{t}{M}\right)g(t,\sigma)\intD t,
\end{align}
that we call windowed time average.
We have that $w(s)\in C^l(\mathbb{R},\mathbb{R})$ is a windowing function if it satisfies
\begin{align}
 w(s)=0,\, s\not\in(0,1),\qquad \int_0^1w(s)\intD s = 1.
\end{align}
Furthermore $l$ denotes the differentiability and we state that $l=-1$ denotes a piecewise continuous function. It is easy to see that the characteristic function of the open interval $(0,1)$, i.e. $\mathds{1}_{(0,1)}$ is a feasible $C^{-1}$ window, which we call Square-window. In this case, the definitions of ${J}_w$ and ${J}_M$ coincide.
Since the windowing function $w$ is independent of $\sigma$, the  windowed time averaged sensitivity is given by
\begin{align}\label{eq_wndAvg2}
   \frac{\intD}{\intD \sigma}{J}_w(\sigma,M) =\frac{1}{M}\int_0^M w\left(\frac{t}{M}\right) \frac{\intD}{\intD \sigma}g(t,\sigma)\intD t.
\end{align} 
In their work \citep{Krakos}, Krakos et al. were able to prove that the convergence of ${J}_w$ and $\frac{\intD}{\intD \sigma}{J}_w $ depends on the differentiability of the windowing function. We recall Krakos' Theorem~\ref{theo1} for the sake of completeness. To this end we define the substitution
\begin{align}
  s(t) = \frac{t}{T},\qquad h( s,\sigma)=g(t,\sigma).
\end{align}
Note that $h( s,\sigma)$ has period length one.
\begin{theorem}\label{theo1}
 The windowed time average ${J}_w$ and the windowed time averaged sensitivity  $\frac{\intD}{\intD \sigma}{J}_w$ computed with a window $w( s)\in C^l(\mathbb{R},\mathbb{R})$ spanning $k = \lceil\frac{M}{T}\rceil$ periods converges to ${J}$ with order
\begin{align}
 \abs{{J}_w(\sigma,M)-{J}(\sigma)}&\leq \norm{h}_\infty  \landauO(k^{-p}), \label{eq_convNonSens}\\
\abs{\frac{\intD}{\intD \sigma}{J}_w(\sigma,M)-\frac{\intD}{\intD \sigma}{J}(\sigma)}&\leq \norm{\partial_\sigma h}_\infty \landauO(k^{-p}) + \norm{\frac{1}{T}\frac{\intD T}{\intD \sigma}}_1 \norm{\partial_s h}_\infty \landauO(k^{-(p-1)}), \label{eq_convSens}
\end{align}
where
\begin{align}
 p =\left \{\begin{array}{@{}ll@{}}
    1,      &  l =-1, \\
    l+1,	& l\geq 0, l\text{ even},\\
    l+2,    & l>0, l\text{ odd}. 
    \end{array}\right.
\end{align}
For $f\in C^\infty(\mathbb{R},\mathbb{R})$, we get an exponential rate of convergence.
\end{theorem}

As a result, we get that the distance between the actual time averaged sensitivity $\frac{\intD}{\intD \sigma} {J}$ and the approximation using a Square-window weighting $\frac{\intD}{\intD \sigma} {J}_M$ is at most in $\abs{\frac{1}{T}\frac{\intD T}{\intD \sigma}}\landauO(1)$. This implies very slow convergence or none at all.
Krakos et al. developed two methods for approximating $ {J}$, namely long-time windowing and short-time windowing. We use the long-time windowing approach which has the advantage that it does not need a priori knowledge of the real period length $T(\sigma)$ and directly utilizes Theorem~\ref{theo1}. The idea of long-time windowing is the following. 
First, one computes the solution of the dynamical system of Eq.~\eqref{eq_ODE}, until it exceeds its transient phase and exhibits periodic behavior. Then, we use a smooth enough window function and compute the windowed time average in Eq.~\eqref{eq_wndAvg1} over a finite time span $M$ starting at $t_{tr}$. We achieve convergence as $M\rightarrow\infty$. 
Similarly, we get the algorithm for the long-time window averaged sensitivity by first calculating the instantaneous sensitivity until the transient phase has passed and then computing the windowed  time averaged sensitivity of the objective function in Eq.~\eqref{eq_wndAvg2}. To compute the instantaneous sensitivity we use the automatic differentiation capabilities of SU2 with both forward and adjoint mode. 

\section{The Optimization Framework}\label{sec_OPT}
The simulation of the compressible URANS equations in the high performance fluid dynamic solver SU2 \cite{SU2_big, towards_SU2, Albring_Sagebaum_Gauger, Albring} is generally done by establishing the method of lines, i.e. one first performs the spatial discretization of convective and viscous fluxes and arrives at the coupled system of ordinary differential equations given by Eq.~\eqref{eq_ODE}. The spatial discretization is implemented using a Finite-Volume scheme on a vertex based median-dual grid with several numerical fluxes available, e.g. JST \cite{JST_original}, Roe \cite{ROE_original}, and AUSM \cite{AUSM}. We solve the resulting ODE system with a dual time stepping method. Utilizing a second-order BDF method for time discretization, we first approximate the solution $U$ to Eq.~\eqref{eq_ODE} by
\begin{align}
R^*(u^n; u^{n-1}, u^{n-2},\sigma) = \frac{3}{2\Delta t}u^n +R(u^n) - \frac{2}{\Delta t}u^{n-1}+\frac{1}{2\Delta t}u^{n-2} = 0, \qquad n=1,\dots,N.
\end{align}
Here $\Delta t$ denotes the length of the physical time step, $n$ is the current (physical) time step, and $N$ the number of (physical) time steps in the simulation. The solution vector of each time step $u^n\in\mathbb{R}^{(d+2)m}$ represents the $d+2$ conservative variables of the URANS equations in $d$ spatial dimensions at $m$ grid points of the computational domain. $u$ is additionally dependent on the variable(s) of a turbulence model, e.g. Spalart-Allmaras \cite{SA_turb}, which is omitted for the sake of simplicity in the notation of this work. The residual $R$ (and therefore $R^*$) is via boundary conditions directly and via the solutions $u^n$ indirectly dependent on the design of the considered structure $\sigma\in R^{n_d}$. Here ${n_d}$ denotes the number of control points of the parametrization of the structure, e.g. an FFD-Box parametrization of an airfoil. 
 Therefore, the extended residual $R^*$ is a mapping given by
\begin{align}
 R^*:\mathbb{R}^{(d+2)m}\times\mathbb{R}^{(d+2)m}\times\mathbb{R}^{(d+2)m}\times\mathbb{R}^{n_d}\longrightarrow\mathbb{R}^{(d+2)m}
\end{align}
The initial condition for the BDF method is given by the free-stream values of the flow-field. The flow at initial time has no physical meaning and needs a transient phase  that is assumed to exceed at time step $n_{tr}$, where it exhibits physically meaningful behavior. We further assume that the system exhibits a limit cycle oscillation after $n_{tr}$ time steps. This assumption is reasonable in many cases, however it should be stressed, that the duration of the transient phase may be dependent on the design parameter $\sigma$ and the value of $n_{tr}$ should be chosen big enough.

A dual time step $\tau$ is used to converge the residual $R^*$ to zero. This can be done by employing an implicit Euler method for the (fictitious time) system of differential equations
\begin{align}\label{eq_ficTimeODE}
 \partial_\tau u^n + R^*(u^n) = 0,
\end{align}
yielding
\begin{align}
 u^n_{p+1}-u^n_p+\Delta\tau R^*(u^n_{p+1};u^{n-1}, u^{n-2}) = 0, \qquad p=1,2,3,\dots
\end{align} 
 Using a linearization of $R^*$ around $U^n_{p}$, given by
 \begin{align}
  u^n_{p+1}-u^n_p+\Delta\tau \left[R^*(u^n_{p}) + \left.\partial_u R^*\right|^n_p\left(u^n_{p+1}-u^n_p\right)\right] = 0, \qquad p=1,2,3,\dots
 \end{align}
 the scheme can be written in the form of a fixed-point iteration
\begin{align}
	G^n: \mathbb{R}^{(d+2)m}\times\mathbb{R}^{(d+2)m}\times\mathbb{R}^{(d+2)m}\times\mathbb{R}^{n_d}\longrightarrow\mathbb{R}^{(d+2)m} 
\end{align}
 that marches $U^n_p$ to a fixed point $U^n$, which can be seen as the steady state solution of the  fictitious time ODE given by Eq.~\eqref{eq_ficTimeODE}. The iterator for each physical time step $n=1,\dots,N$ reads as
\begin{align}
	u^n_{p+1}=G^n \left(u^n_p,u^{n-1},u^{n-2},\sigma\right), \qquad p=1,2,3,\dots\, .
\end{align}
 $U^{n-1}$ and $U^{n-2}$ denote the fixed-points of the previous physical time steps.
 To set up the discrete optimal control problem, we approximate the continuous ouput in Eq.~\eqref{eq_wndObj} by a midpoint rule.
 \begin{align}
      {J}_w(\sigma,M) =\frac{1}{M}\int_0^M w\left(\frac{t}{M}\right)g(t,\sigma)\intD t \approx \frac{1}{N-n_{tr}}\sum_{n=n_{tr}}^N w\left(\frac{n-n_{tr}}{N-n_{tr}}\right) g(U^n(\sigma),\sigma)
 \end{align}
We can now formulate the discrete optimal control problem
\begin{subequations}\label{eq_OCP}
\begin{align}
	\min_{\sigma\in X_{ad}}&\qquad	  {J}_w = \frac{1}{N-n_{tr}}\sum_{n=n_{tr}}^N w\left(\frac{n-n_{tr}}{N-n_{tr}}\right) g(U^n(\sigma),\sigma)\label{eq_wndObj} \\
	\text{subject to}&\qquad U^n(\sigma) = G^n\left(U^n,U^{n-1},U^{n-2},\sigma\right), \qquad n = 1,\dots,N, \label{eq_state}
\end{align}
\end{subequations}
where we use  the midpoint rule for a discrete approximation to $ {J}_w$, which we also denote by $J_w$. The time span to take the average is given by $M=N-n_{tr}$. Note that taking the limit $M\rightarrow\infty$ in the long-time windowing approach is equivalent to taking the limit $N\rightarrow\infty$ in Eq.~\eqref{eq_wndObj}. The physical time step, where the system has passed the transient phase, is denoted by $n_{tr}$. 
The set of admissible designs is given by $X_{ad}\subset\mathbb{R}^{n_d}$. 

The idea to solve the optimization problem is to evaluate the objective function with a flow solution corresponding to an initial design $\sigma_0$. We then compute a descent direction of the objective function, that is given by the total derivative of $ {J}_w$ with respect to $\sigma$, $-\frac{\intD}{\intD \sigma} {J}_w$, that is then used in an update step using e.g. a quasi Newton method or a line search method.  The implementation of  the SU2 solver for discrete optimization problems with PDE constraints is done via the adjoint approach, that we want to adjust with respect to the weighted objective function $ {J}_w$. Here, the idea is to avoid the computation of computationally expensive terms which are marked red in Eq.~\eqref{eq_sens}, which appear when applying the chain rule to $\frac{\intD}{\intD \sigma} {J}_w$. 
\begin{align}\label{eq_sens}
	\frac{\intD}{\intD \sigma} {J}_w\left([U^1,\dots,U^N](\sigma),\sigma\right) = \partial_\sigma {J}_w\left([U^1,\dots,U^N],\sigma\right) + \sum_{n=1}^N \partial_{U^n} {J}_w \textcolor{red}{\frac{\intD U^n}{\intD \sigma}}
\end{align}
We therefore introduce the Lagrangian corresponding to the optimal control problem given by Eq.~\eqref{eq_OCP}. This is a function 
\begin{align}
	L:\mathbb{R}^{N(d+2)m}\times\mathbb{R}^{N(d+2)m}\times\mathbb{R}^{n_d}\longrightarrow\mathbb{R}
\end{align} 
defined as
\begin{align}\label{eq_lag}
\begin{aligned}
&L\left(\left[U^1,\dots,U^N\right]^T,\left[\overline{U}^1,\dots,\overline{U}^N\right]^T,\sigma\right)\\
	 &= \frac{1}{N-n_{tr}}\sum_{n=n_{tr}}^N w\left(\frac{n-n_{tr}}{N-n_{tr}}\right) J(U^n(\sigma),\sigma) -
	  \sum_{n=0}^N \left[ \left(\overline{U}^n\right)^T\left( U^n(\sigma) - G^n\left(U^n,U^{n-1},U^{n-2},\sigma\right)\right)\right], 
\end{aligned}
\end{align}
where  $\overline{U}^n$ denotes the adjoint state vector at time step $n$.
 We obtain the following KKT system under the assumption of differentiability of $L$ with repect to $\sigma$, $U^n$, and $\overline{U}^n$ for all $n=1,\dots,N$, that is given by the construction of the fixed-point iteration, 
\begin{align}
 \partial_{\overline{U}^n}L &\overset{!}{=} 0, \qquad n= 1,\dots,N,\qquad \text{state equations}\\
 \partial_{U^n}L &\overset{!}{=} 0, \qquad n=1,\dots,N,\qquad \text{adjoint equations}\label{eq_adj}\\
 \partial_\sigma L  &\overset{!}{=} 0, \qquad \qquad\qquad\quad\qquad\,\, \text{design equation.}
\end{align}
The state equations are already given by Eq.~\eqref{eq_state}. The discrete adjoint equations, that are used to efficiently compute the total derivative of the objective function $ {J}_w$ with respect to $\sigma$, are obtained by computing the partial derivatives $\partial_{U^n}L$ for all $n=1,\dots,N$.  The adjoint equations of the KKT system, given by Eq.~\eqref{eq_adj}, yield the expression
\begin{align}\label{eq_discAd}
	\overline{U}^n = \left(\partial_{U^n} G^{n}\right)^T \overline{U}^n  + \left(\partial_{U^n} G^{n+1}\right)^T \overline{U}^{n+1} +   \left(\partial_{U^n} G^{n+2}\right)^T \overline{U}^{n+2}
	+ \mathds{1}_{\menge{n\geq n_{tr}}}\frac{1}{N-n_{tr}}\left(w\left(\frac{n-n_{tr}}{N-n_{tr}}\right)\partial_{U^n}J(U^{n})\right)^T,
\end{align}
for all $n=N,\dots,1$. The characteristic function $\mathds{1}_{\menge{n\geq n_{tr}}}$ indicates, that the seeding of the objective function is only performed from time step $n_{tr}$ on.
Remark that since the window function $w$ is only dependent on the current time step $n$, the adjustment to traditional adjoint state equations \cite{Albring} is marginal.  This eases the integration into existing high performance solvers like SU2. Analogously to the direct equations we may write Eq.~\eqref{eq_discAd} as a set of fixed-point iterations $N_U^n$
\begin{align}\label{eq_adjIter}
	\overline{U}_{p+1}^{n} &= H_U^n\left(\overline{U}_p^n,\overline{U}^{n+1}, \overline{U}^{n+2}, \sigma\right)  \\
	&= \left(\partial_{U^n} G^{n}\right)^T \overline{U}_p^n  + \left(\partial_{U^n} G^{n+1}\right)^T \overline{U}^{n+1} +   \left(\partial_{U^n} G^{n+2}\right)^T \overline{U}^{n+2}
	+ \mathds{1}_{\menge{n\geq n_{tr}}}\frac{1}{N-n_{tr}}\left(w\left(\frac{n-n_{tr}}{N-n_{tr}}\right)\partial_{U^n}J(U^{n})\right)^T \nonumber.
\end{align}
For all $ p=1,2,3,\dots $. This adjoint fixed-point iteration has to converge for each physical time step $n$, starting from final time $N$ and marching backwards in time like the continuous adjoint equation formulated in \cite{Nadarajah2007OptimumSD}. The adjoint states $\overline{U}^{n+1}$ and $\overline{U}^{n+2}$ denote the converged adjoint fixed-points of the previous time steps $n+2$ and $n+1$.
By the Banach fixed-point Theorem we have convergence of this fixed-point iteration, if it is a contraction. The norm of the derivative of the iterator $H^n_U$ with respect to its argument $\overline{U}^n$ can be written as 
\begin{align}
	\norm{\partial_{\overline{U}^{n}}N_U^n} = 
	\norm{\left(\partial_{U^n} G^{n}\right)^T}. 
\end{align}
Convergence of the adjoint iterator therefore depends only on the convergence of the primal iterator.
We have $\norm{\left(\partial_{U^n} G^{n}\right)^T} < 1$ in a suitable norm, if the direct (pseudo time) iteration is near convergence to a steady state solution.
Finally, the design equation is given by the partial derivative of the Lagrangian with respect to the design $\sigma$,
\begin{align}\label{eq_design_eq}
	\partial_\sigma L = \sum_{n=0}^N\left[\mathds{1}_{\menge{n\geq n_{tr}}}\frac{1}{N-n_{tr}}w\left(\frac{n-n_{tr}}{N-n_{tr}}\right)\partial_\sigma J(U^{n}) + \left(\overline{U}^n\right)^T\partial_\sigma G^n \right].
\end{align}

\section{Periodic detached flow around the NACA0012 airfoil}\label{subsec_TestCaseA}
In the following, we present numerical results of an aerodynamic shape optimization performed in the test case of the NACA0012 airfoil. All computations are done on the RHRK high performance computing center in Kaiserslautern. The presented methods are published in the open source CFD suite SU2, version 7.0.1 "Blackbird".

The test case at hand is a subsonic unsteady turbulent flow around the NACA0012 airfoil. The mesh is a quadrilateral O-grid that wraps around the NACA0012 airfoil. It consists of $27125$ elements, of which $217$ are wall boundary elements and $217$  are farfield boundary elements. We use the URANS solver of the SU2 CFD suite, where we employ the Jameson-Schmidt-Turkel (JST) \citep{JST_original} scheme for the Navier-Stokes fluxes, whereas the turbulent viscosity is calculated using the Spalart-Allmaras turbulence model \cite{SA_turb} with a scalar Upwind scheme.
A FFD-Box parametrization is used to model the airfoil surface with a total of 242 control points, that act as design variables for the sensitivity computations. In particular, for the comparison between the windows, a single FFD design variable is chosen.
The used flow configuration is displayed in the following.
\begin{itemize}
 \item Mach Number: $0.3$
 \item Angle of Attack: $17.0$° (Degrees)
 \item Freestream Temperature: $293.0 K$
 \item Reynolds Number: $10^3$, $10^6$
 \item Reynolds Length: $1.0$ (Length of the airfoil)
\end{itemize} 
The simulation is computed with time step $\Delta t = 0.0005s$ for several final times $t_{f}$. Naturally, the discrete objective function is computed with the same time step.

\begin{figure}
\begin{minipage}{.5\textwidth}
  \centering
  \includegraphics[width=.9\linewidth]{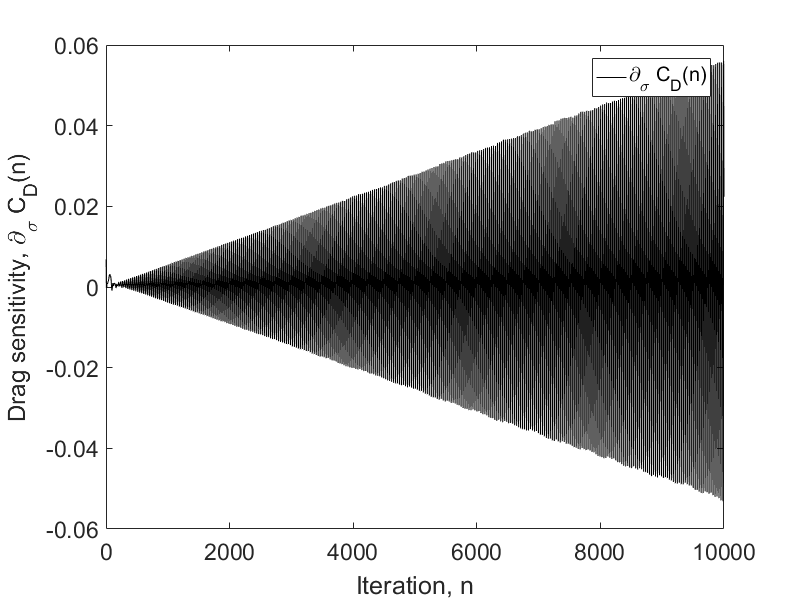}{\\ (a) linear scale, $Re=10^3$.}
\end{minipage}
\begin{minipage}{.5\textwidth}
  \centering
  \includegraphics[width=.9\linewidth]{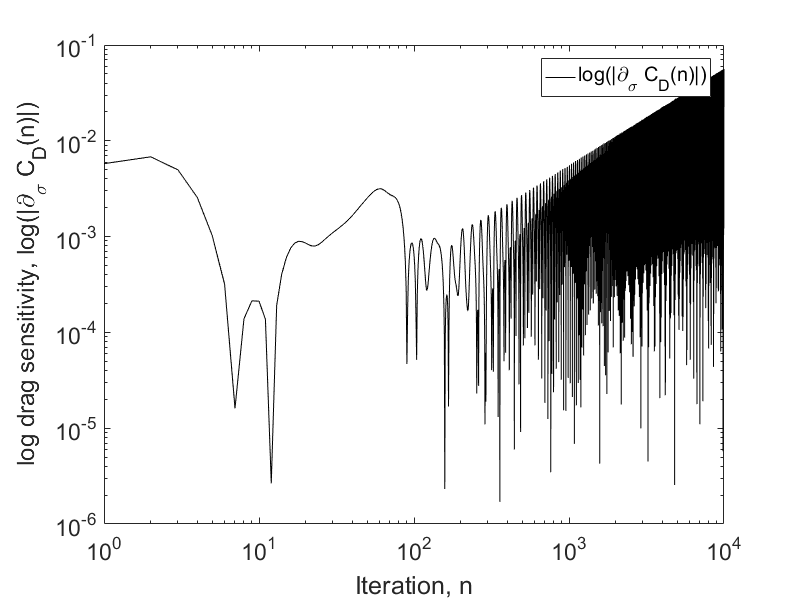}{\\ (b) logarithmic scale, $Re=10^3$.}
\end{minipage}
\begin{minipage}{.5\textwidth}
  \centering
  \includegraphics[width=.9\linewidth]{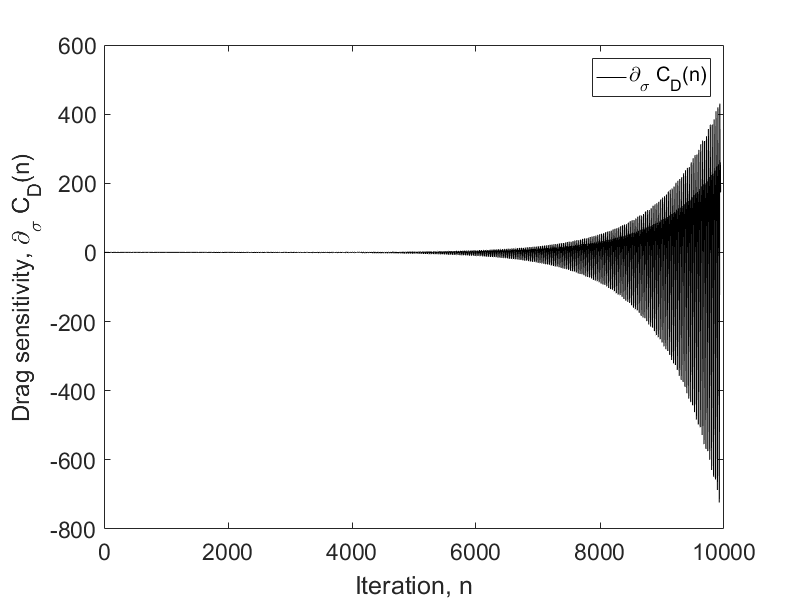}{\\ (c) linear scale, $Re=10^6$.}
\end{minipage}
\begin{minipage}{.5\textwidth}
  \centering
  \includegraphics[width=.9\linewidth]{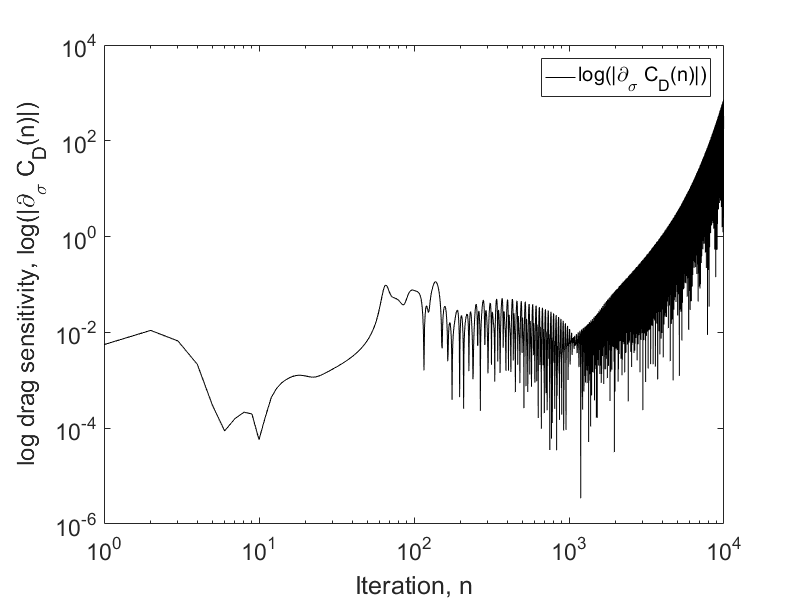}{\\ (d) logarithmic scale, $Re=10^6$.}
\end{minipage}
  \captionof{figure}{Long-time behavior of $\frac{\intD}{\intD \sigma} C_D(n)$  over $n$; higher Reynolds numbers, i.e. more turbulent flows lead to exponential growths in the period amplitude of the sensitivity}
  \label{fig_RE}
\end{figure}

\begin{figure}
\begin{minipage}[t]{.5\textwidth}
  \centering
  \includegraphics[width=.9\linewidth]{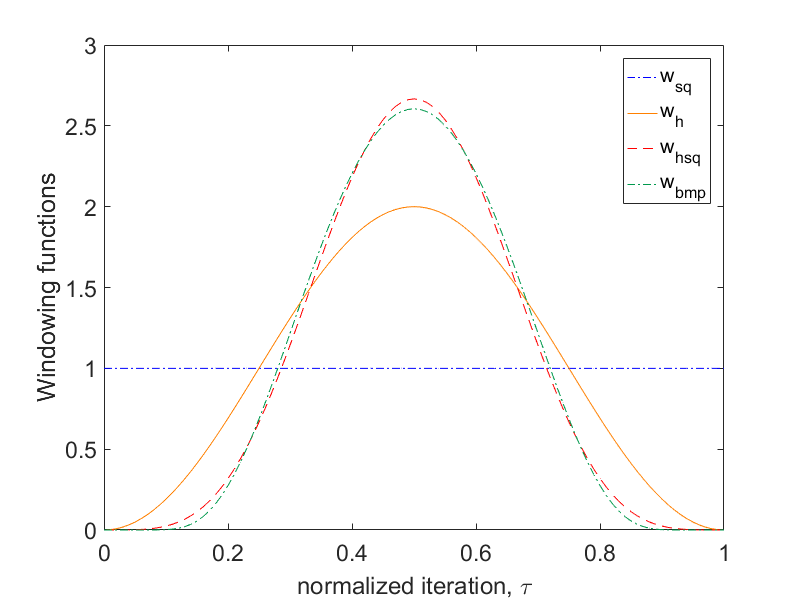}
  \captionof{figure}{Square- , Hann-, Hann-Square-, and Bump-window over $s$}
  \label{fig_wndFcts}
\end{minipage}\qquad
\begin{minipage}[t]{.5\textwidth}
  \centering
  \includegraphics[width=.9\linewidth]{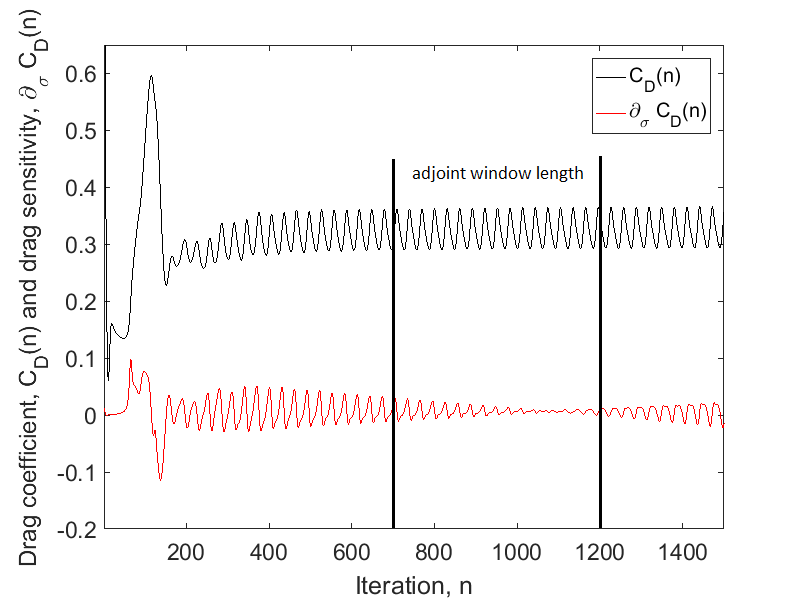}
  \captionof{figure}{$C_D(n)$ and $\frac{\intD}{\intD \sigma} C_D(n)$  over $n$ with $Re=10^6$; the window for the adjoint run spans iteration $500$ to $1200$}
  \label{fig_CL}
\end{minipage}
\end{figure}

The sensitivity $\frac{\intD C_D}{\intD \sigma}$ of the drag coefficient $C_D$ with respect to the design variables is computed with help of the built in automatic differentiation (AD) support of SU2 based on the software package CoDiPack \cite{Sagebaum}. We use two modes of AD. First, the primal (forward) mode, that represents a tangent vector evaluation. The computational effort scales with the number of design variables. The second, more important mode is the adjoint mode of AD, that is independent of the number of design variables. To achieve an efficient sensitivity calculation, advanced AD techniques such as checkpointing are used in combination with the above described Lagrangian framework. In Subsection~\ref{sec_valWND}, the forward AD mode is used, however in Subsection~\ref{sec_valAD} we show consistency of both approaches in our test case and thus the results of Subsection~\ref{sec_valWND} are also valid for the more practical adjoint mode.
\begin{figure}
\centering
\begin{minipage}{0.33\textwidth}
  \centering
  \includegraphics[width=.95\linewidth]{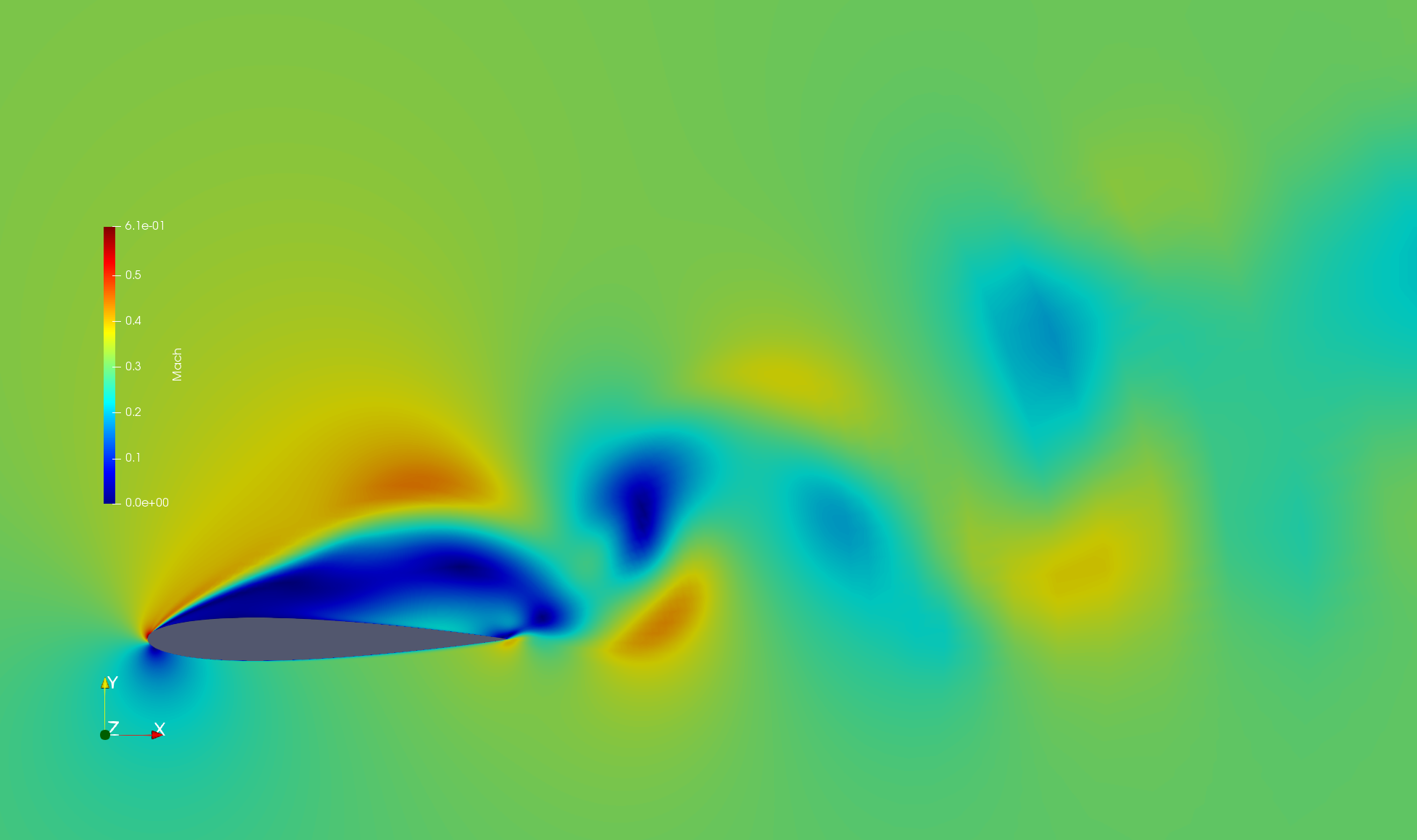}{\\(a) $t_1 = 1/3T(\sigma)$}
\end{minipage}
\begin{minipage}{0.33\textwidth}
  \centering
  \includegraphics[width=.95\linewidth]{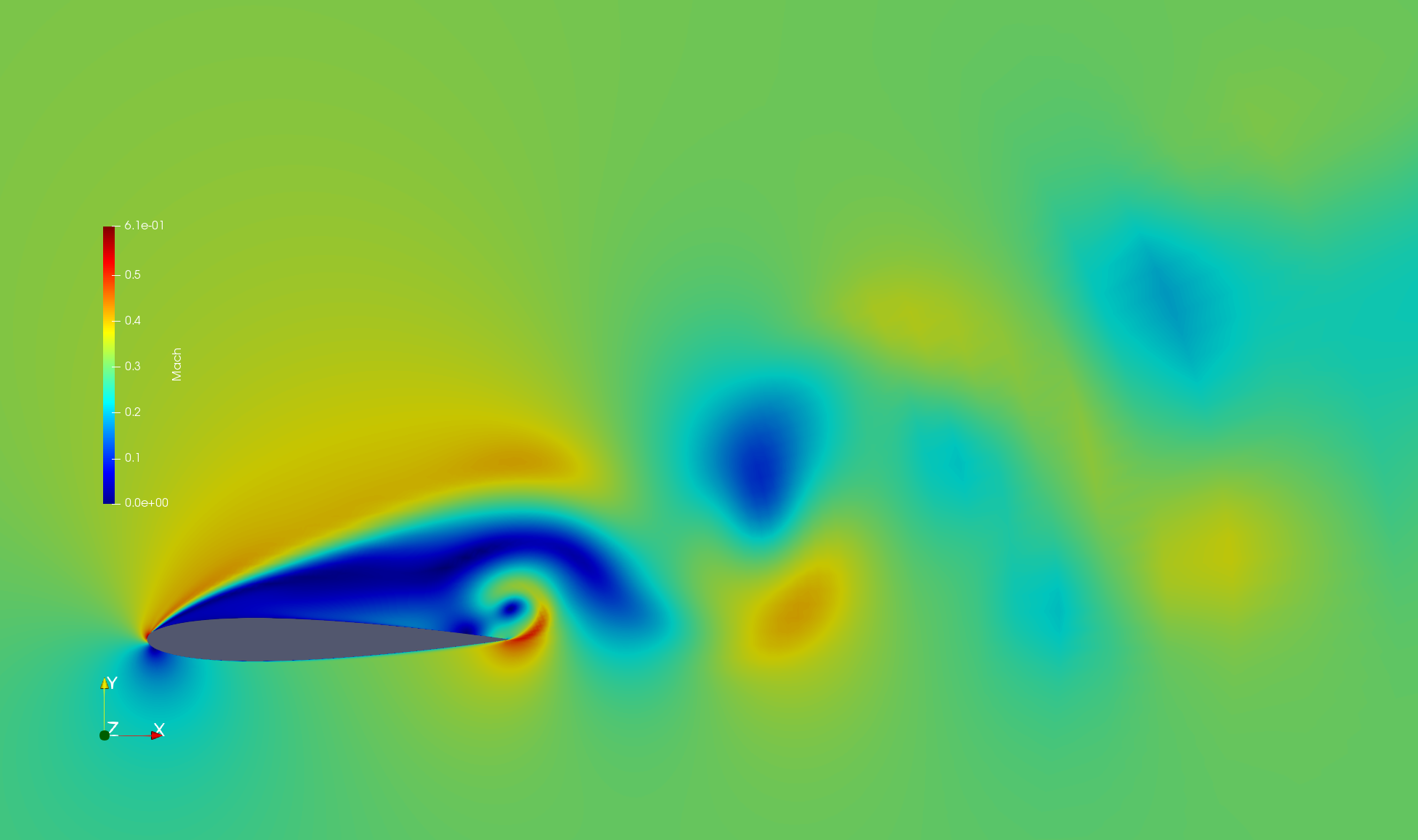}{\\(b) $t_2 = 2/3T(\sigma)$}
\end{minipage}
\begin{minipage}{0.33\textwidth}
  \centering
  \includegraphics[width=.95\linewidth]{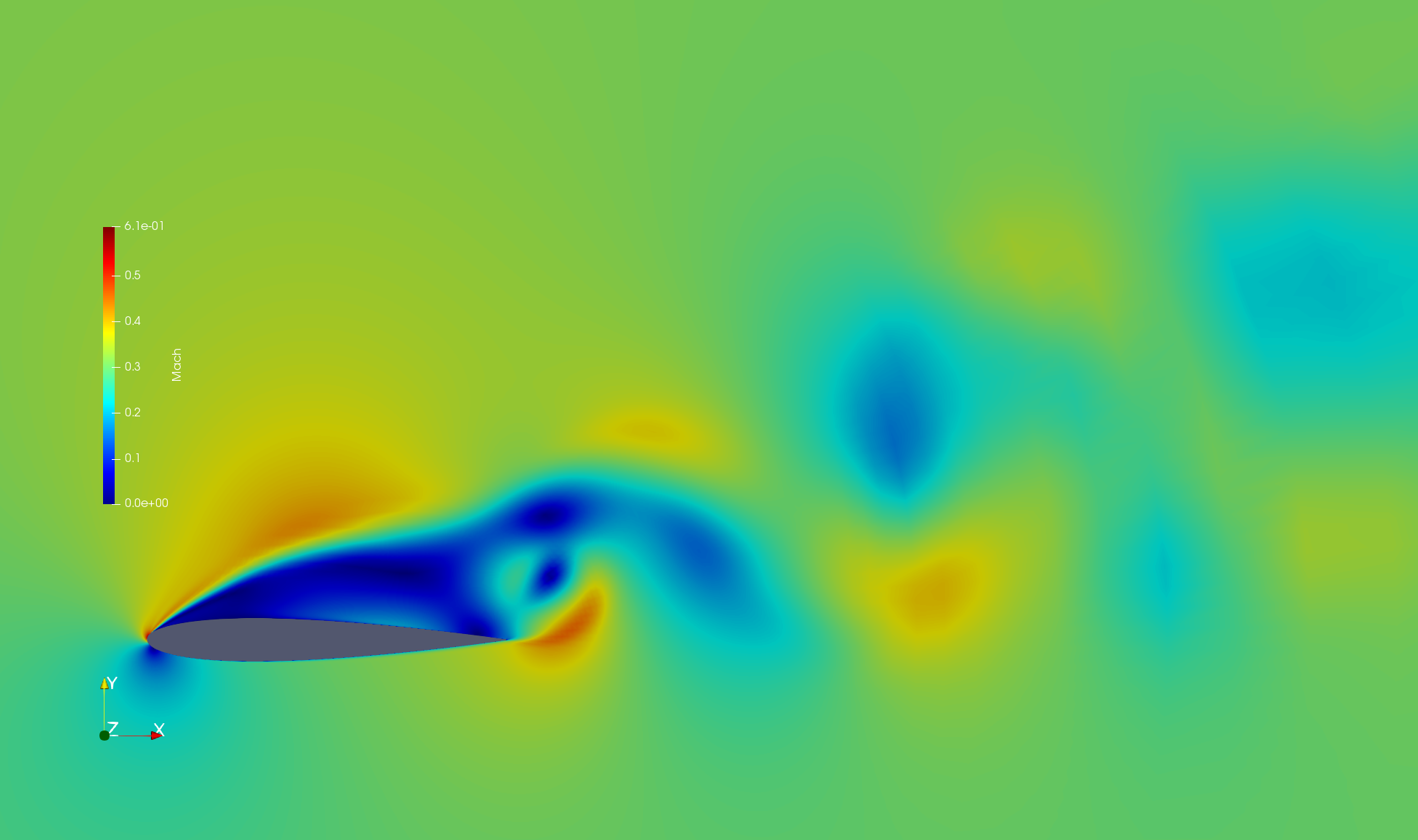}{\\(c) $t_3 = T(\sigma)$}
\end{minipage}%
 \captionof{figure}{Velocities of a flow around NACA0012 airfoil with $Re = 10^6$, $17$ degrees AoA at different time points of a period. The period length is $T(\sigma)\approx 31 \Delta t$.}
 \label{fig_nacaFlow}
\end{figure}
The angle of attack of $17$° results in a detached  primal flow, as can be seen in Fig.~\ref{fig_nacaFlow}. In the wake of the airfoil, a vortex street can be seen that illustrates the oscillating character of the flow, which is exhibited after the transient phase. We can see the transient phase in Fig.~\ref{fig_CL} as well as the periodic behavior of $C_D$ in its limit cycle oscillation. It should be noted, that the sensitivity $\frac{\intD}{\intD \sigma} C_D$ exhibits periodic behavior in the same time frame with the same period length as $C_D$. However, the amplitude of its limit cycle oscillation changes in time. Therefore, it is worthwhile to consider first the long-time behavior of the drag sensitivity with different Reynolds numbers, which is displayed in Fig.~\ref{fig_RE}.

We compare the flow configuration above with two choices of Reynolds numbers, $Re = 10^3$ and $Re = 10^6$. Both flows are detached and exhibit periodic behavior after some transient phase. Considering the drag sensitivities, we can see in Fig.~\ref{fig_RE}a and  Fig.~\ref{fig_RE}b that the amplitude of $\frac{\intD}{\intD \sigma} C_D$ grows linearly in the case of $Re=10^3$. However, in the case of $Re=10^6$ the growth in the amplitude is exponential, as can be seen in  in Fig.~\ref{fig_RE}a and  Fig.~\ref{fig_RE}b.

Taking Theorem~\ref{theo1} into account, we expect good convergence behavior for the case $Re=10^3$ and worse behavior for the case $Re=10^6$. This is due to the fact, that the upper bound for the difference $\abs{\frac{\intD}{\intD \sigma} {J}_w(\sigma,M)-\frac{\intD}{\intD \sigma} {J}(\sigma)}$ scales with $\norm{\partial_s C_D}_\infty$ and $\norm{\frac{1}{T}\frac{\intD T}{\intD \sigma}}_1\norm{\partial_\sigma C_D}_\infty$. If $\norm{\frac{\intD}{\intD \sigma} C_D}_\infty$ grows exponentially, by the chain rule $\norm{\partial_s C_D}_\infty$ or $\norm{\frac{1}{T}\frac{\intD T}{\intD \sigma}}_1\norm{\partial_\sigma C_D}_\infty$ must grow exponentially as well, which can affect the convergence behavior of the windowed time averaged sensitivity. However, in our test case, the exponential growth of the amplitude dominates only after some time, see Fig.~\ref{fig_RE}b, so it is still worth considering the windowed sensitivities.

\subsection{Validation of the windowing approach}\label{sec_valWND}
In the following we validate the embedding of the windowing approach in the sensitivity analysis that is necessary for aerodynamic shape optimization. We apply several windows to the time averaged lift coefficient $C_D$
\begin{align}\label{eq_wndAvg}
  {J}_w(\sigma,t_{f}-t_{tr}) =\frac{1}{t_{f}-t_{tr}}\int_{t_{tr}}^{t_{f}} w\left(\frac{t-t_{tr}}{t_{f}-t_{tr}}\right)C_D(t,\sigma)\intD t,
\end{align}
where $t_{f}$ is an arbitrarily chosen final time and $t_{tr}$ is a starting time at which the flow has exceeded its transient phase. Analogously, we consider the discrete windowed time averaged lift coefficient
\begin{align}\label{eq_discWndAvg}
  {J}_w(\sigma,N-n_{tr}) =\frac{1}{N-n_{tr}}\sum_{n=n_{tr}}^{N} w\left(\frac{n-n_{tr}}{N-n_{tr}}\right)C_D(n,\sigma).
\end{align}
Here $N=\frac{t_{f}}{\Delta t}$ is the final time step, where $\Delta t$ is the time step size,  and $n_{tr}$ is the time step corresponding to time $t_{tr}$.
We compare several windowing functions, which were suggested by Krakos et al. \cite{Krakos}, with different orders of differentiability, namely
\begin{align}
 \text{Square-window:}& \qquad w_{sq}(s) = \mathds{1}_{(0,1)}(s), \qquad &s\in(0,1),\\
 \text{Hann window:}& \qquad w_{h}(s) = 1-\cos(2\pi s), \qquad &s\in(0,1),\\
 \text{Hann-Square-window:}& \qquad w_{hsq}(s) = \frac{2}{3}\left(1-\cos(2\pi s)\right)^2, \qquad &s\in(0,1),\\
 \text{Bump-window:}& \qquad w_{bmp}( s) =\frac{1}{A}\exp\left(\frac{-1}{ s- s^2}\right), \qquad & s\in(0,1), 
\end{align}
where
\begin{align}
 A=\int_0^1 \exp\left(\frac{-1}{s-s^2}\right)\intD s.
\end{align}
The different windowing functions can be seen in Fig.~\ref{fig_wndFcts}. Considering the order of differentiability, we get the following orders of convergence $p$ for the windowed time average and order $p_s$ for the windowed time averaged sensitivity.
\begin{itemize}
 \item $w_{sq}\in C^0 \Rightarrow p=1$ and $p_s=0$.
 \item $w_{h}\in C^1 \Rightarrow p=3$ and $p_s=2$.
 \item $w_{hsq}\in C^3 \Rightarrow p=5$ and $p_s=4$.
 \item $w_{bmp}\in C^\infty \Rightarrow p=\infty$ and $p_s=\infty$, i.e. exponential rate of convergence.
\end{itemize} 
It is reasonable to start the windowed time averaging at time $t_{tr}=500$ where the transient phase of the system has already passed and the limit cycle oscillator has been reached. One period of $C_D$ consists of approximately $31$ iterations in both configurations. 

\begin{figure}
\begin{minipage}{.5\textwidth}
  \centering
  \includegraphics[width=.9\linewidth]{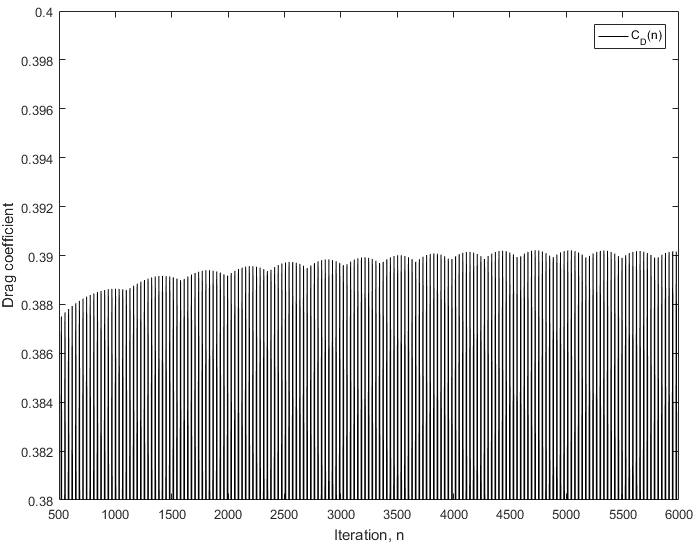}{\\ (a) $Re = 10^3$.}
\end{minipage}
\begin{minipage}{.5\textwidth}
  \centering
  \includegraphics[width=.9\linewidth]{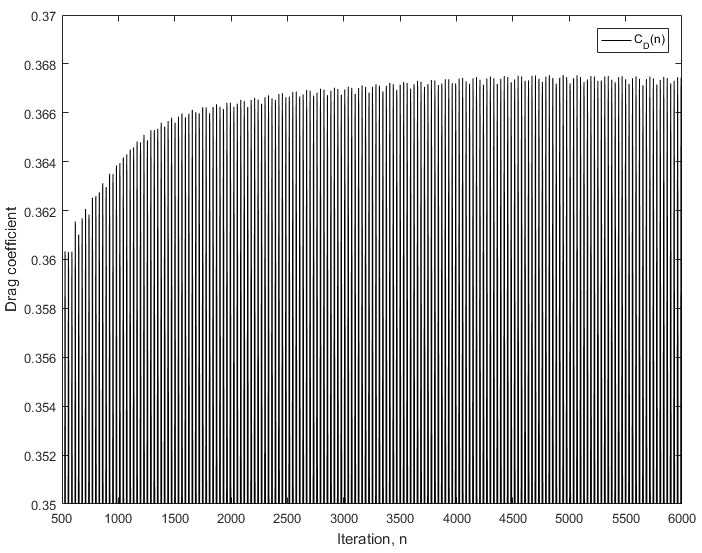}{\\ (b) $Re = 10^6$.}
\end{minipage}
\caption{ $C_D(n)$ long-time behavior with different Reynolds numbers. Note the slight shift upwards of the period mean in both cases, that stops at approximately $n=4500$. }
\label{fig_shift}
\end{figure}

\begin{figure}
\begin{minipage}{.5\textwidth}
  \centering
  \includegraphics[width=.9\linewidth]{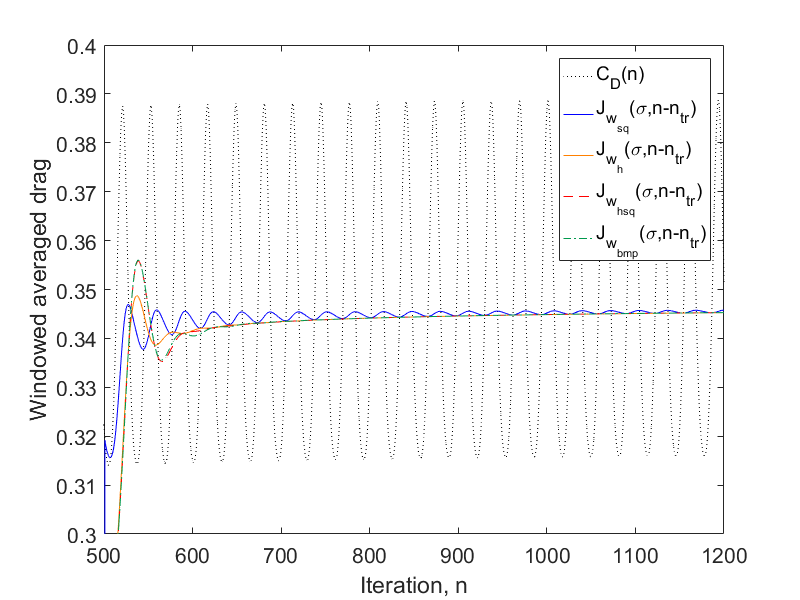}{\\ (a) $C_D(n)$ and $ {J}_w(\sigma,n-n_{tr})$ up to time step $1200$.}
\end{minipage}
\begin{minipage}{.5\textwidth}
  \centering
  \includegraphics[width=.9\linewidth]{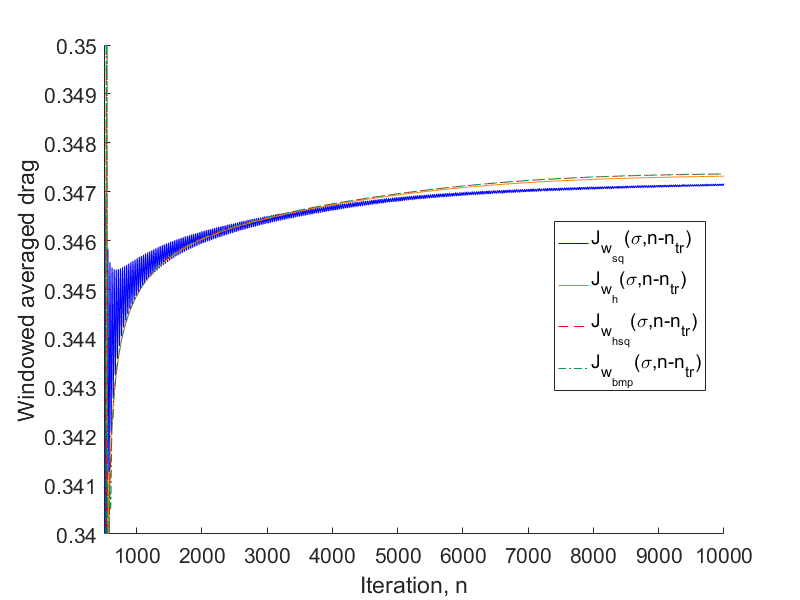}{\\ (b) $ {J}_w(\sigma,n-n_{tr})$ up to time step $10000$.}
\end{minipage}
\caption{ $C_D(n)$ and $ {J}_w(\sigma,n-n_{tr})$ with different windows over $n$, $n_{tr}=500$, $Re = 10^3$.}
\label{fig_wndCDRe3}
\end{figure}

\begin{figure}
\begin{minipage}{.5\textwidth}
  \centering
  \includegraphics[width=.9\linewidth]{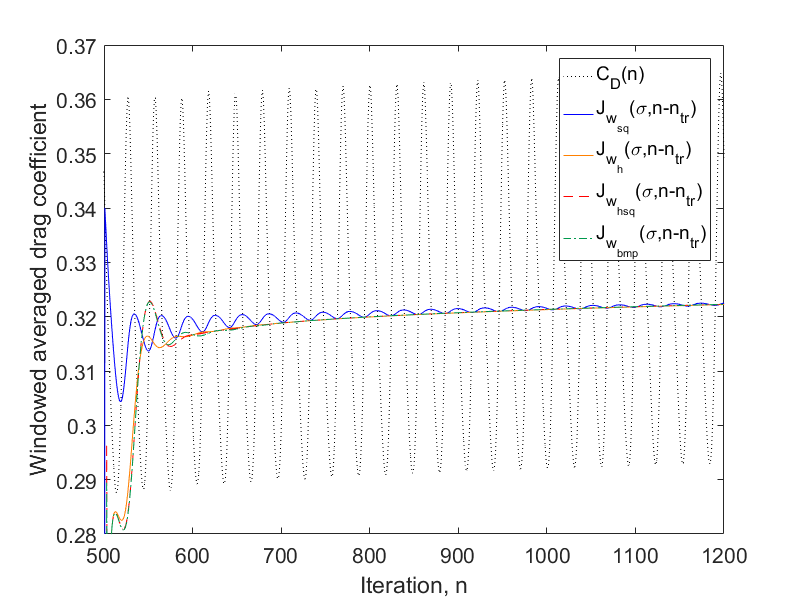}{\\ (a) $C_D(n)$ and $ {J}_w(\sigma,n-n_{tr})$ up to time step $1200$.}
\end{minipage}
\begin{minipage}{.5\textwidth}
  \centering
  \includegraphics[width=.9\linewidth]{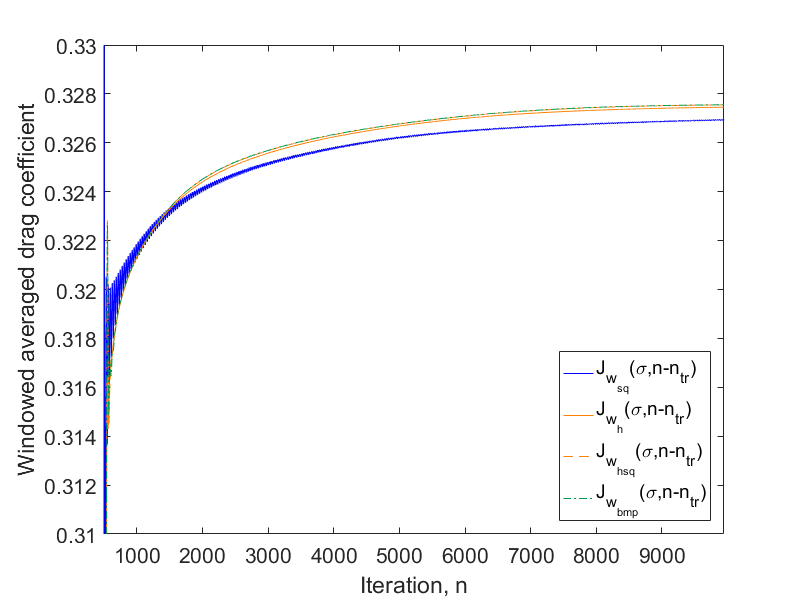}{\\ (b) $ {J}_w(\sigma,n-n_{tr})$ up to time step $10000$.}
\end{minipage}
\caption{ $C_D(n)$ and $ {J}_w(\sigma,n-n_{tr})$ with different windows over $n$, $n_{tr}=500$, $Re = 10^6$.}
\label{fig_wndCL}
\end{figure}

We start the analysis of our results with the windowed drag coefficient for both Reynolds cases $Re=10^3$ and $Re=10^6$, since the results are similar.
Figure~\ref{fig_wndCDRe3}a respectively Fig.~\ref{fig_wndCL}a shows $C_D(n)$ over time step $n$ with above described flow configuration as well as the windowed time averages over end time step $n$. Hence, the time span to average $C_D$ is $M=n-500$. We can observe, that all high order windows converge within a short amount of periods. The Bump-window needs approximately $k=5$ periods to level of, the Hann windowed average oscillates for a some more periods before converging to the same value as the Bump-windowed average. We can see that the amplitude of the Square-window also reduces, but at a much slower rate than the higher order windows.  Note that the windowed time averages for final times $n\leq 531$, i.e. where the time span is less than one period, are not meaningful. 

 Now we consider the long-time behavior of different windowed time averages for final times $n=10^4$ which is equal to $k=300$ periods in Fig.~\ref{fig_wndCDRe3}b respectively Fig.~\ref{fig_wndCL}b. We can observe several effects. First, none of the higher order windows shows any oscillatory behavior, whereas the Square-window still oscillates for very high iteration numbers. Second, we can see an asymptotic convergence of all windows as the iteration number increases. Both effects are expected, since Theorem~\ref{theo1} postulates convergence in $\landauO(k^1)$ for the Square-window, where $k$ is the number of periods passed within the averaged time span.
Furthermore, we can see in Fig.~\ref{fig_shift}a respectively Fig.~\ref{fig_shift}b that the mean value of the period shifts slightly upwards. In  Fig.~\ref{fig_wndCDRe3}b respectively Fig.~\ref{fig_wndCL}b, we can see that all windowed time averages follow this trend. Indeed we have observed in Fig.~\ref{fig_shift} that the shift in the mean value stops at $n=4500$. We can see that higher order windows have less problems adapting to this trend in their long-time behavior, since they weight function values at (normalized) times near to $s=0$ and $s =1$ much less than the Square-window, as can be seen in Fig.~\ref{fig_wndFcts}. More generally speaking, the decreased weighting of values near to $s=0$ and $s =1$ leads to a decreased sensitivity with respect to the behavior of the instantaneous output in these regions. As a result, higher order windowed averages output useful values although in this example there exists an additional trend, i.e. the function is not exactly periodic.

\begin{figure}
\begin{minipage}{.5\textwidth}
  \centering
  \includegraphics[width=.9\linewidth]{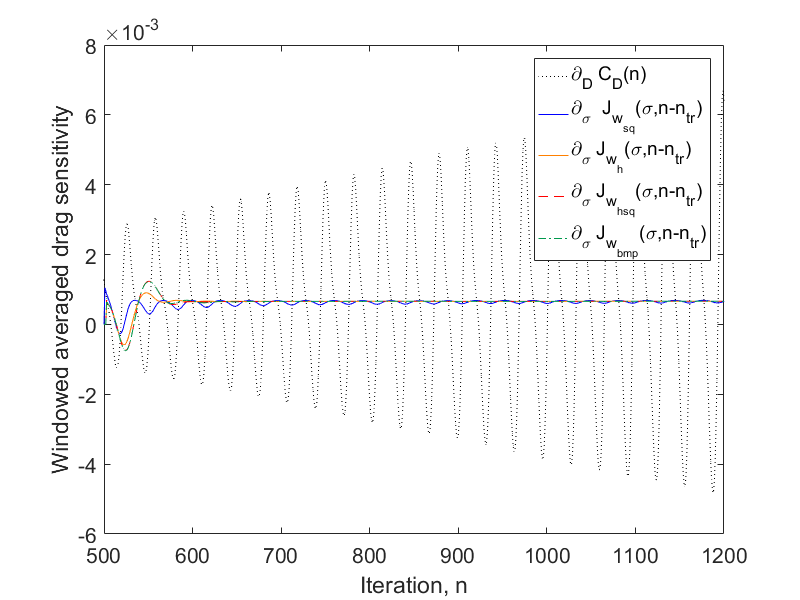}{\\(a)  $\frac{\intD}{\intD \sigma}C_D(n)$ and $\frac{\intD}{\intD \sigma} {J}_w(\sigma,n-n_{tr})$ up to time step $1200$.}
\end{minipage}
\begin{minipage}{.5\textwidth}
  \centering
  \includegraphics[width=.9\linewidth]{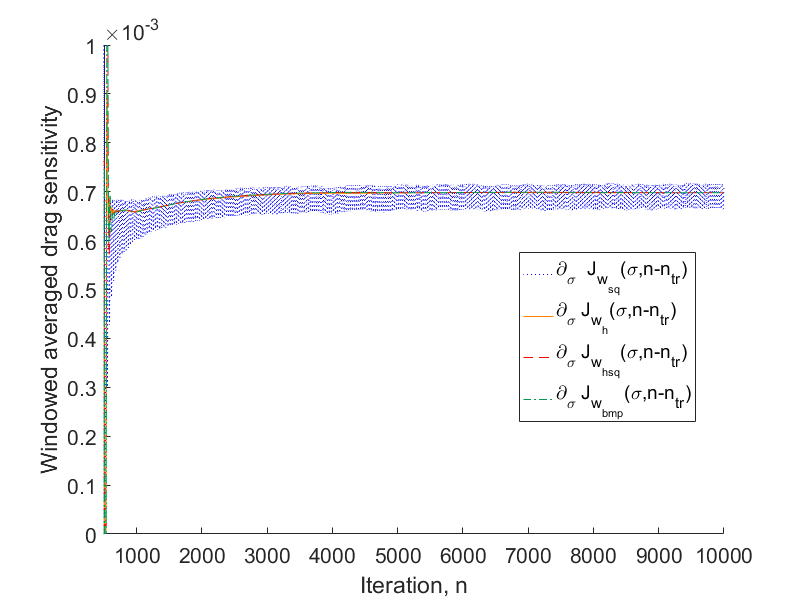}{\\(b)  $\frac{\intD}{\intD \sigma} {J}_w(\sigma,n-n_{tr})$ up to time step $5000$.}
\end{minipage}
 \caption{$\frac{\intD}{\intD \sigma}C_D(n)$ and $\frac{\intD}{\intD \sigma} {J}_w(\sigma,n-n_{tr})$ with different windows over $n$, $n_{tr}=500$, $Re = 10^3$.}
 \label{fig_sensOSCRe3}
\end{figure}

\begin{figure}
\begin{minipage}{.5\textwidth}
  \centering
  \includegraphics[width=.9\linewidth]{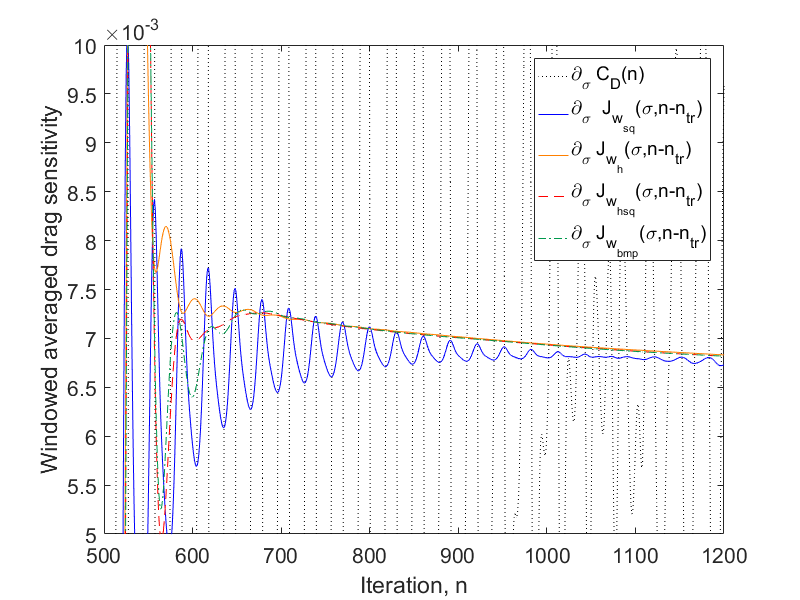}{\\(a)  $\frac{\intD}{\intD \sigma}C_D(n)$ and $\frac{\intD}{\intD \sigma} {J}_w(\sigma,n-n_{tr})$ up to time step $1200$.}
\end{minipage}
\begin{minipage}{.5\textwidth}
  \centering
  \includegraphics[width=.9\linewidth]{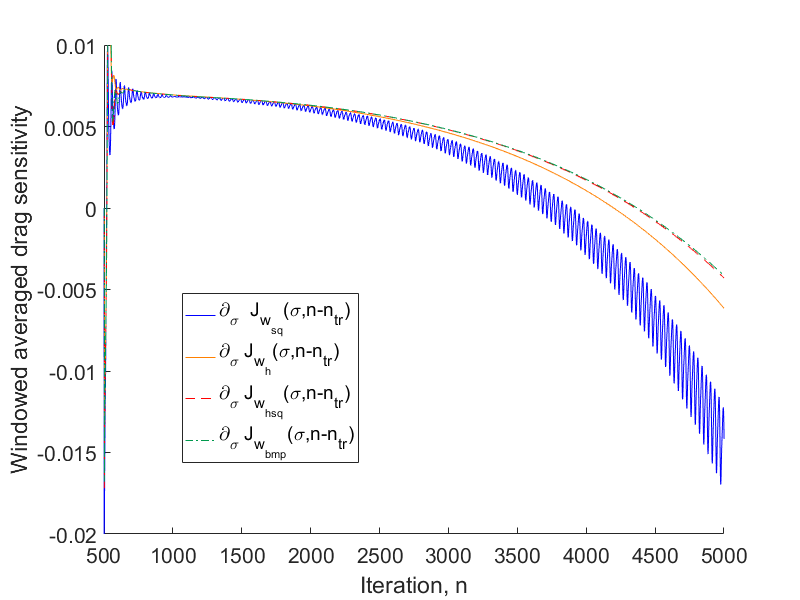}{\\(b)  $\frac{\intD}{\intD \sigma} {J}_w(\sigma,n-n_{tr})$ up to time step $5000$.}
\end{minipage}
 \caption{$\frac{\intD}{\intD \sigma}C_D(n)$ and $\frac{\intD}{\intD \sigma} {J}_w(\sigma,n-n_{tr})$ with different windows over $n$, $n_{tr}=500$, $Re = 10^6$.}
 \label{fig_sensOSC}
\end{figure}

In analogy to Fig.~\ref{fig_wndCL}, the sensitivity of $C_D$ with respect to the design variable $\sigma$ as well as the corresponding windowed time averages are displayed in Fig.~\ref{fig_sensOSC} with $Re=10^6$ and in Fig.~\ref{fig_sensOSCRe3} with $Re=10^3$. 
Let us first analyze the results of the lower Reynolds number case, i.e. the less chaotic system. The results completely validate Theorem~\ref{theo1}, since the higher order windowed time averaged sensitivities quickly converge, see Fig.~\ref{fig_sensOSCRe3}b, and the Square-window leads to an oscillatory time average with a non-decreasing amplitude. Theorem~\ref{theo1} postulates a convergence rate for the Square-windowed time average in $\landauO(1)$. This is in line with the findings of Krakos et. al \cite{Krakos}, who analysed convergence behavior of the windows.

Now let us consider the more chaotic flow with $Re=10^6$.
Considering Fig.~\ref{fig_sensOSC}a, we can observe again that the higher order windowed averaged sensitivities do not show any oscillatory behavior after the first few periods.
The Square-window first reduces its amplitude, but this is solely due to the fact, that the amplitude of the instantaneous sensitivity reduces up to $n=1050$. We can see in Fig.~\ref{fig_sensOSC}b,  that the amplitude of the Square-window increases again with increasing iteration number $n$. This reflects the findings of Theorem~\ref{theo1}, that postulates no convergence for the Square-windowed time averaged sensitivity.
Now we consider the long-time behavior of the higher order windows in Fig.~\ref{fig_sensOSC}b. We can observe a divergent behavior for all windows as time increases. The divergence rate is slow for the higher order windows, but higher for the Square-window. This is a result of the above discussed feature of higher order windows to be less sensitive to values of the instantaneous output near to $s=0$ and $s=1$.
 The windowed time averaged sensitivities assume smaller and smaller values, but considering the corresponding (instantaneous) sensitivities in Fig.~\ref{fig_RE}c, we can see that the limes inferior of $\frac{\intD}{\intD \sigma} C_D$, grows much faster to $-\infty$ than the limes superior of $\frac{\intD}{\intD \sigma} C_D$ grows to $\infty$. Thus, the effective mean value of a period also shifts downwards, which is reflected by the divergent behavior of the windowed time averaged sensitivities. However, it should be stressed, that Theorem~\ref{theo1} does not give clear convergence bounds in this case as discussed in Subsection~\eqref{subsec_TestCaseA}, since the exponential growth of the amplitude of the sensitivity implies exponential growth of $\norm{\frac{\intD}{\intD \sigma} C_D}_\infty$.

Also interesting is the time frame between $n=3500$ and $n=4000$ iterations, where the Square-windowed time averaged sensitivity changes its sign due to its oscillation. This can lead to massive errors in the sensitivity computation. In the context of aerodynamic shape optimization, this may result in wrong signs or magnitude of the components of the sensitivity vector, i.e. a non-diminishing error in the descent direction of the design process. In such cases, a higher order window that does not oscillate yields much more meaningful results, that can be used in an shape optimization process.


\subsection{Validation of the adjoint solver}\label{sec_valAD}
 
The adjoint flow solver of SU2 is validated using the NACA0012 airfoil with the same flow conditions as above once with $Re=10^3$ and once with $Re=10^6$. To compare adjoint and primal computation, we choose a fixed final time step $N=1200$ and a fixed starting time $n_{tr} = 500$. We average therefore over $700$ time steps which corresponds to approximately $22.5$ periods.  Table~\ref{tab_ad_vs_dd} shows the values of the sensitivity with respect to one design variable, computed with forward and adjoint mode, as well as their relative difference in percent (\%). By forward mode we refer to the primal (tangent) computation of the derivatives using the forward AD mode. By adjoint mode, we refer to the Lagrangian method derived in Section~\ref{sec_OPT} in combination with the adjoint AD mode.
 The adjoint derivative calculation is performed using an iterative scheme, that utilizes the fixed-point structure of the flow solver, i.e. Eq.~\eqref{eq_adjIter} and Eq.~\eqref{eq_design_eq}. Therefore, consistency of adjoint and primal sensitivities depend on the convergence of the adjoint iteration in Eq.~\eqref{eq_adjIter}.
 
 We can see that the relative difference of both computation methods is a maximum of $1.88$\% in case of the Hann-window at $Re=10^6$, but mostly below $0.5$\%. The Hann-Square- window at $Re=10^6$ displays the lowest relative difference with $0.0318$\%. As a conclusion, we can say that both methods to compute the sensitivities are consistent. This enables the efficient adjoint method for the shape optimization procedure.
 
 \begin{table}
\centering
    \begin{tabular}{| l | l | l | l | l |}
    \hline
    Reynolds number & Windowing Function & Adjoint Mode  & Forward Mode & Relative Difference, \%\\ \hline
   $10^6$& Square           & 0.00770318 	& 0.00769755 & 0.0731 \\ \hline
   $10^6$& Hann             & 0.00819720 	& 0.00835195 & 1.8878 \\ \hline
   $10^6$& Hann-Square      & 0.00836323 	& 0.00836589 & 0.0318 \\ \hline
   $10^6$& Bump             & 0.00839826 	& 0.00836385 & 0.4097 \\ \hline
   $10^3$& Square 			& 0.00065344	& 0.00065595 & 0.3841 \\ \hline
   $10^3$& Hann   			& 0.00066978	& 0.00066094 & 1.3198 \\ \hline
   $10^3$& Hann-Square 		& 0.00067026	& 0.00066131 & 1.3353 \\ \hline
   $10^3$& Bump 			& 0.00067049 	& 0.00066155 & 1.3344 \\ \hline
    \end{tabular}
    \caption{Sensitivity of the drag coefficient $C_D$ with respect to a single FFD-Box design variable, computed with different differentiation techniques and different windows.}
    \label{tab_ad_vs_dd}
\end{table}

\subsection{Comparison of surface sensitivities using different windows}\label{subsec_surfsens}
In the following, we consider the surface sensitivities of $C_D$ in the setting of the NACA0012 airfoil. We choose surface sensitivities over sensitivities with respect to the design variables, since they give more intuitive insights into how the wing is deformed in the design process. In this test case, there are $217$ surface points on the airfoil. 
We compare surface sensitivities calculated with the Bump- and Square-window with different windowing time spans. The Square-window acts as a benchmark and represents the case where no specific windowing function is chosen and therefore the traditional approach.

As discussed in Subsection~\ref{sec_valWND} and shown in Fig.~\ref{fig_sensOSC}a, the Square-windowed time average first reduces and then increases its amplitude, which distorts the convergence behavior of the Square-window. Therefore, we choose $n_{tr}=1503$ instead of $n_{tr}=500$. 
We average over $700$ and $709$ iterations up to final iteration $n_{f1}=2203$ respectively $n_{f2}=2212$. The difference of $9$ iterations is approximately  $29$\% of a period length.
\begin{figure}
\centering
\begin{minipage}{\textwidth}
  \centering
  \includegraphics[width=\linewidth]{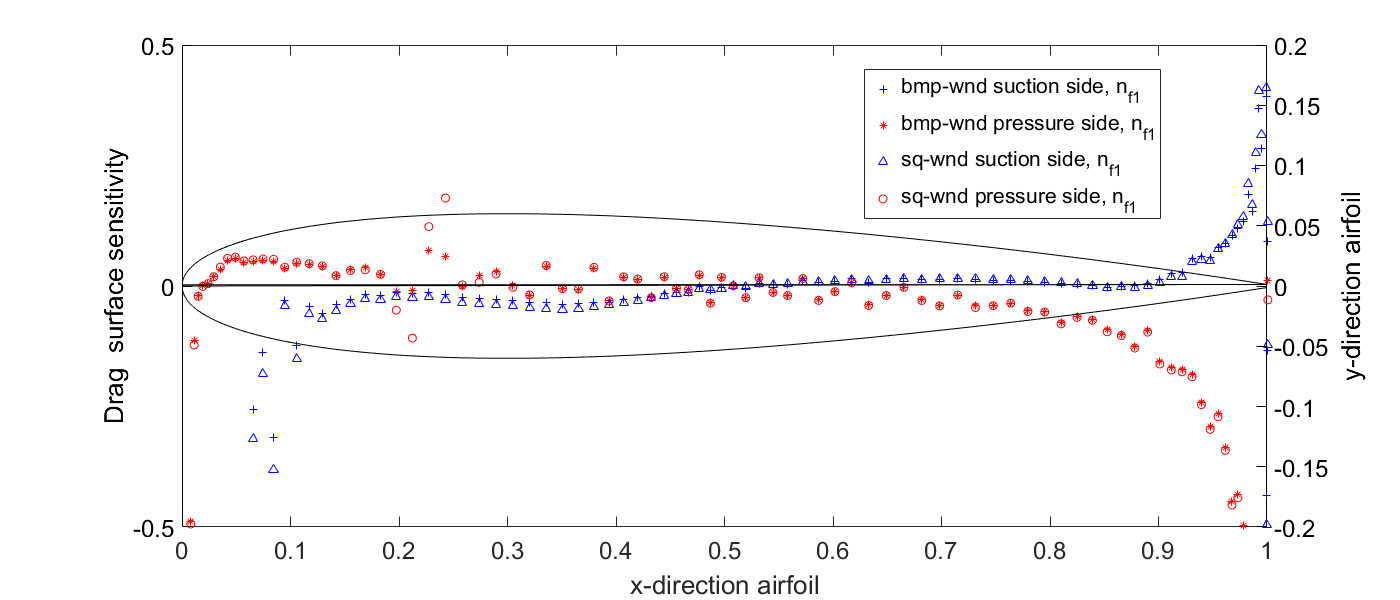}{\\(a) Bump- and Square-window at $n_{f1}$.}
\end{minipage}
\begin{minipage}{\textwidth}
  \centering
  \includegraphics[width=\linewidth]{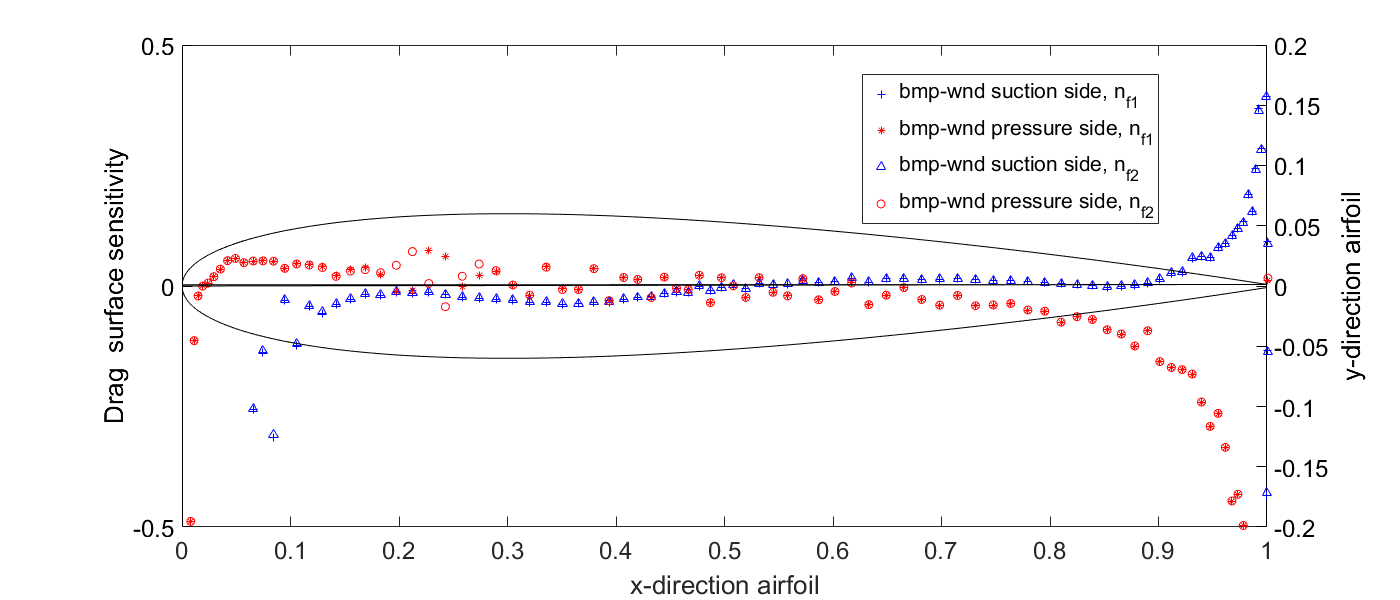}{\\(b) Bump-window at $n_{f1}$ and $n_{f2}$.}
\end{minipage}
\begin{minipage}{\textwidth}
  \centering
  \includegraphics[width=\linewidth]{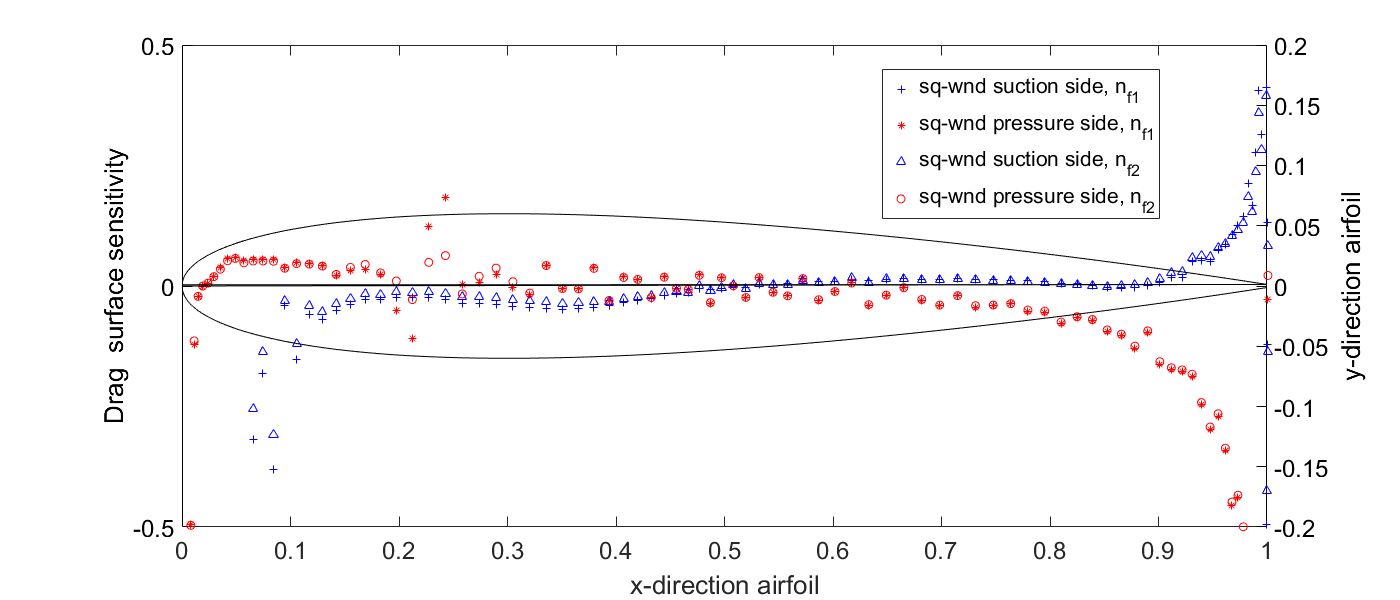}{\\(c)  Square-window at $n_{f1}$ and $n_{f2}$.}
\end{minipage}
\caption{Drag surface sensitivities at $Re=10^6$, $n_{tr}=1503$ with different end times $n_{f1}=2203$, $n_{f2}=2212$.}
\label{fig_surfSens}
\end{figure}
Figure~\ref{fig_surfSens} shows Bump- and Square-windowed time averaged surface sensitivities with different final times $n_{f}$. The blue symbols represent the sensitivities on the suction side, whereas the red symbols represent those on the pressure side of the airfoil. The sensitivities are plotted on their respective x-axis locations on the airfoil, i.e. for every point on the x-axis, there are 4 sensitivities plotted. For reference, the NACA0012 airfoil is plotted on the right hand side y-axis.

In Fig.~\ref{fig_surfSens}a, we can see Bump- and Square-windowed time averaged surface sensitivities at final time $n_{f1}$. In both methods, the sensitivities at the suction side near the leading edge are the highest, followed by the sensitivities on both sides near the trailing edge. This is expected, since the flow detaches at the leading edge and forms big vortices at the trailing edge, as can be seen in Fig.~\ref{fig_nacaFlow}. Note that the suction side sensitivities near the leading edge are too big in magnitude to be displayed in Fig.~\ref{fig_surfSens}. 
The absolute difference between both methods is the highest near the leading edge at suction side followed by a region at $20$\% airfoil length, where the Square-window values show a deviation from the Bump-window values. Except for this region, the relative difference of the sensitivities ranges from $5-20$\%. 

Next, we consider Fig.~\ref{fig_surfSens}b, where the Bump-window values are shown for different final times $n_{f1}$ and $n_{f2}$. We can see that even at the high sensitivity regions at the trailing and leading edge, the Bump-windowed time averaged sensitivities barely differ. This shows that the Bump-window is quite robust with respect to the choice of time span to average. 
Contrary to this, we can see in Fig.~\ref{fig_surfSens}c, that the Square-windowed time averaged sensitivities change clearly over the span of $n_{f2}-n_{f1}=9$ iterations. 
In the highly sensitive area at $20$\% airfoil length, the changes in the Square-window values are big compared to the changes in the Bump-window values in Fig.~\ref{fig_surfSens}b.
The relative difference  of the sensitivity vectors in Euclidean norm computed with Bump- and Square-window at $n_{f1}$ is $9.2$\%, which displays a significant change in the descent direction of the optimization procedure. We can even observe a sign change at $5$ surface points, when comparing the two windows. The relative difference between the Bump- windowed sensitivity vectors at $n_{f1}$ and $n_{f2}$ is $0.84$\% in comparison to the relative difference between the Square-windowed sensitivity vectors at both times, that comes to $9.1$\%. 

On the one hand, typically one does not know the exact period length of the objective function and therefore is forced to choose an arbitrary time span to average the windowed sensitivity; on the other hand, a relatively small change in the chosen time span leads to a significant error in the Square-windowed time averaged sensitivity.
This further reinforces the thesis that higher order windows, such as the Bump-window,
leads to a more robust sensitivity calculation and therefore to  increased robustness in the shape optimization with respect to the choice of the averaging time span.

\subsection{Comparison of optimization results using different windows}

In this subsection, we evaluate the results of an aerodynamic shape optimization run performed by  {SU2} with different choices of windowing functions as objective function regularizers. We consider the optimization problem
\begin{equation}\label{eq_ASP_EX}
\begin{aligned}
\min_{\sigma\in X_{ad}}\quad &{J}_w(\sigma,N-n_{tr}) =\frac{1}{N-n_{tr}}\sum_{n=n_{tr}}^{N} w\left(\frac{n-n_{tr}}{N-n_{tr}}\right)C_D(n,\sigma)\\
\text{s.t.}\quad &G(U^n,\sigma)=U^n , \qquad  \forall n = 1,\dots,N     \\
				 &{C}_w(\sigma,N-n_{tr}) =\frac{1}{N-n_{tr}}\sum_{n=n_{tr}}^{N} w\left(\frac{n-n_{tr}}{N-n_{tr}}\right)C_L(n,\sigma)\geq 0.96,\\
				 &X_{ad} = [-0.05,0.05]^{242}\in\mathbb{R}^{242}.
\end{aligned}
\end{equation}
Here, $n_d=242$ and $X_{ad}$ is cuboid in $R^n_d$. $G(u^n,\sigma)$ is the fixed point form of the the URANS solver, which we eliminate formally by the Lagrangian in Eq.~\eqref{eq_lag}. We input the obtained sensitivities to the {SLSQP} optimization algorithm by Kraft \cite{SLSQP}, which then computes a new design for the next design iteration.

\begin{figure}
\centering
\begin{minipage}{0.8\textwidth}
  \centering
  \includegraphics[width=\linewidth]{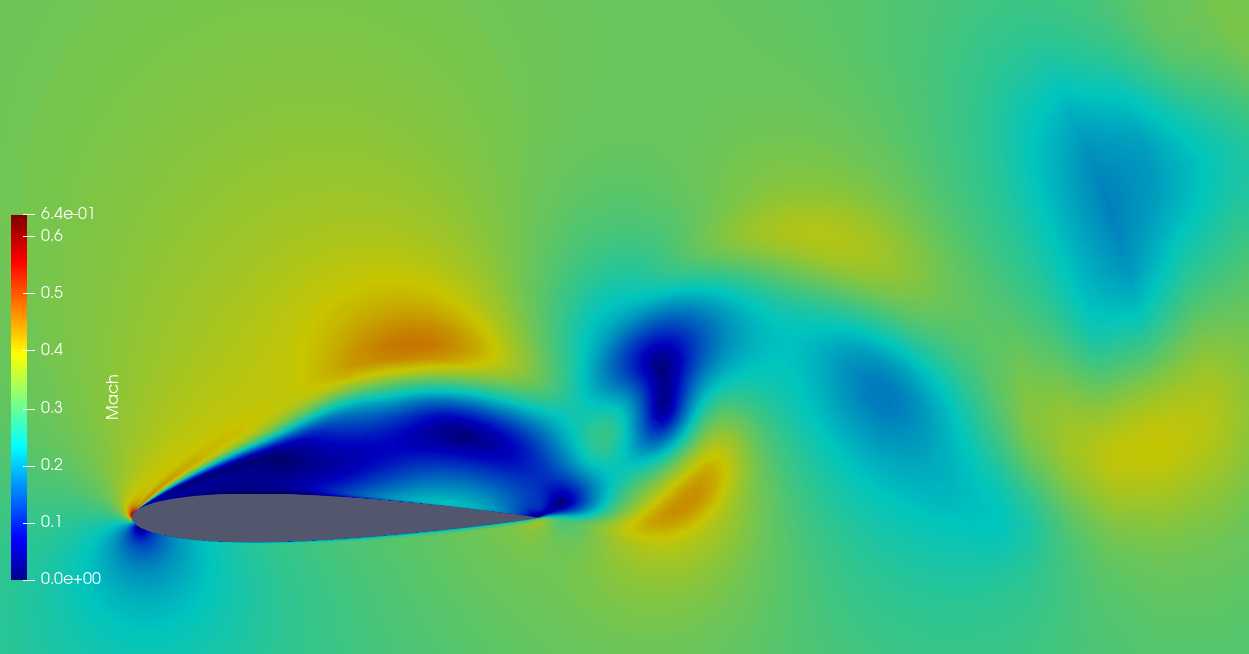}{\\(a) Flow around the baseline design given by the NACA0012-airfoil.}
\end{minipage}
\begin{minipage}{0.8\textwidth}
  \centering
  \includegraphics[width=\linewidth]{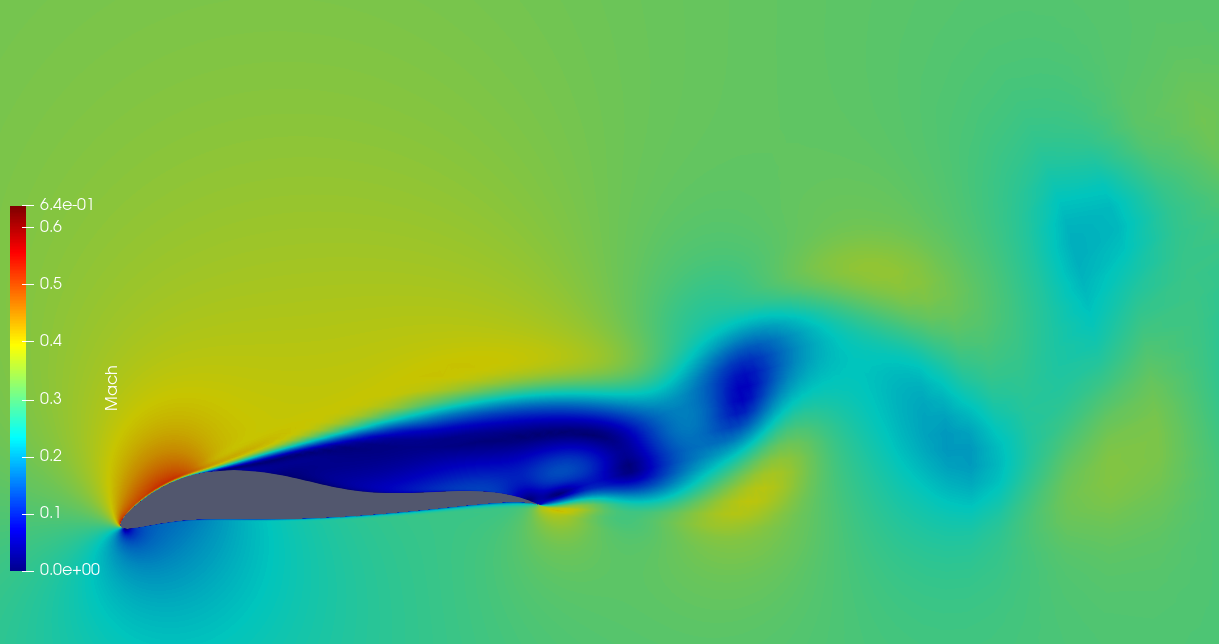}{\\(b) Flow around optimized design using Square-windowing.}
\end{minipage}
\begin{minipage}{0.8\textwidth}
  \centering
  \includegraphics[width=\linewidth]{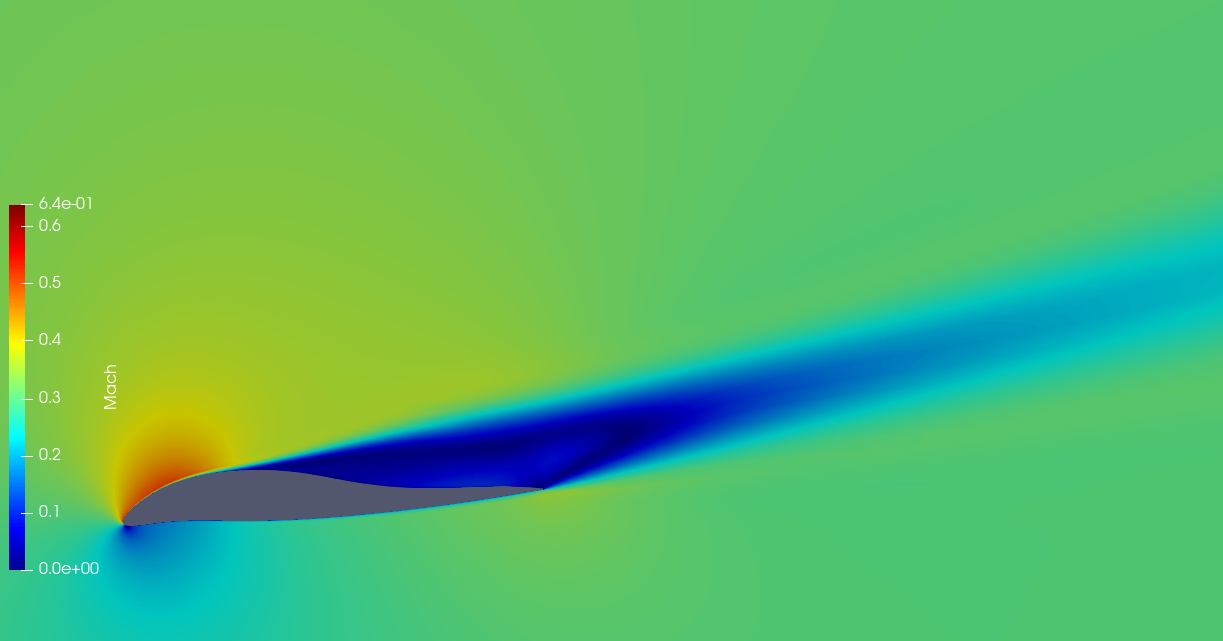}{\\(c)  Flow around optimized design using Hann-Square-windowing.}
\end{minipage}
\caption{Comparison of baseline design and optimized designs using different windows. Displayed are the Mach number of each control volume at time step $n=1500$.}
\label{fig_optflow}
\end{figure}

\subsubsection*{Test Case $Re=10^6$}
 As inequality constraint we use the windowed time averaged lift coefficient, that should be greater or equal to $0.96$, which is the approximate value of the windowed time averaged lift  for the baseline configuration. We start the windowed time average at iteration $n_{tr} = 1500$ and average up to $N=2200$. This corresponds to the configuration in the last subsection. The flow configuration is  described in the beginning of this chapter with Reynolds number $Re=10^6$. We multiply the computed gradients with a relaxation factor of $0.1$, such that the step sizes proposed by the optimizer are not too big for the admissible set $X_{ad}$. 

\begin{figure}
\begin{minipage}{.5\textwidth}
  \centering
  \includegraphics[width=.9\linewidth]{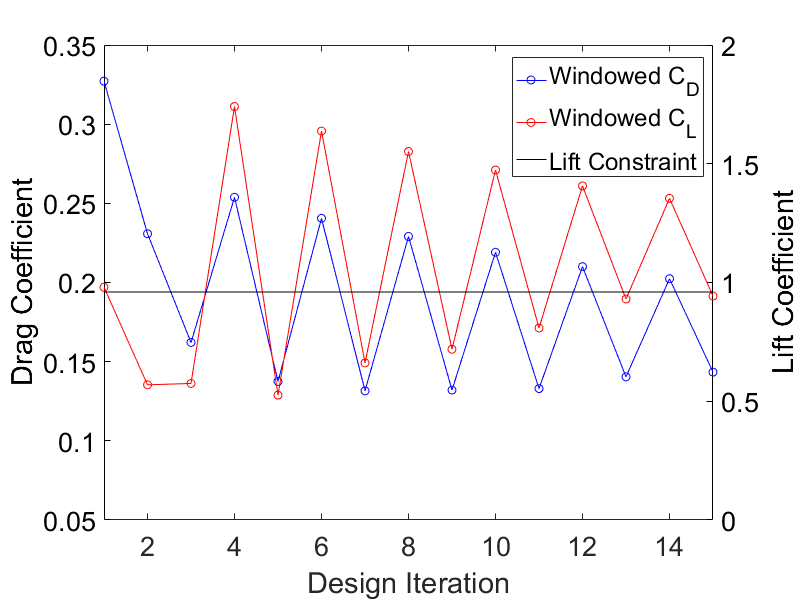}{\\(a) Square-windowing.}
\end{minipage}
\begin{minipage}{.5\textwidth}
  \centering
  \includegraphics[width=.9\linewidth]{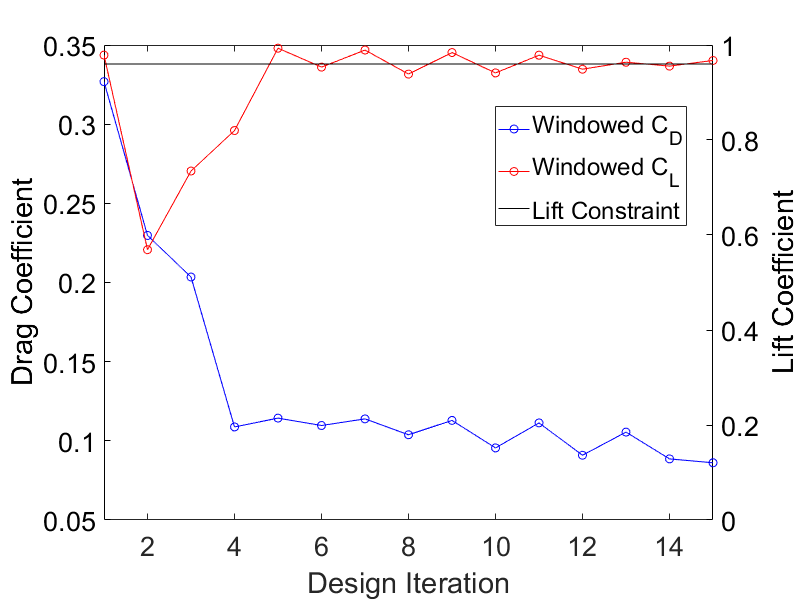}{\\(b) Hann-windowing.}
\end{minipage}
\begin{minipage}{.5\textwidth}
  \centering
  \includegraphics[width=.9\linewidth]{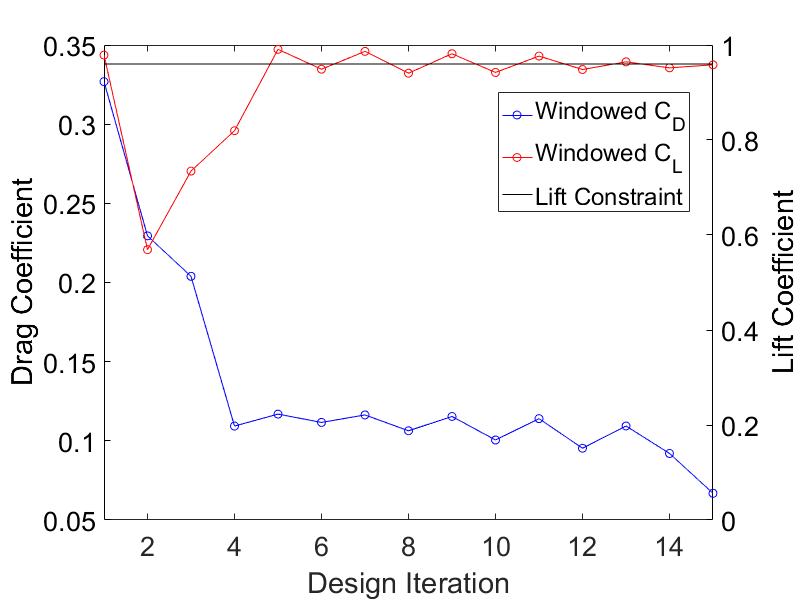}{\\(c) Hann-Square-windowing.}
\end{minipage}
\begin{minipage}{.5\textwidth}
  \centering
  \includegraphics[width=.9\linewidth]{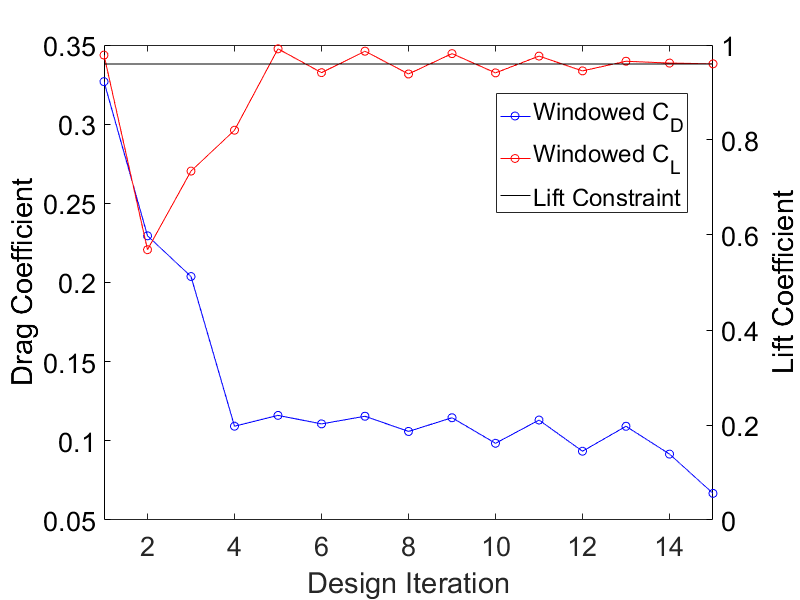}{\\(d)  Bump-windowing.}
\end{minipage}
 \caption{Shape optimization performed by SLSQP with different windows used. Both - feasible and infeasible - design steps are displayed. $Re=10^6$.}
 \label{fig_dsgn_res}
\end{figure}

We compare the optimization results of optimization procedures computed with the different windowing functions introduced above. Figure~\ref{fig_optflow} compares the flow around the baseline geometry described by the NACA0012 airfoil to the results of the shape optimization process, where once the traditional Square-window  and once the Hann-Square-window is used. All flow field snapshots are taken at time step $n=1500$ i.e. past the transient phase of the flow around the base geometry. Whereas the flow optimized using the Square-window, depicted in Fig.~\ref{fig_optflow}b, is unsteady and periodic, the design optimized with the Hann-Square-window exhibits a stationary flow, as seen in Fig.~\ref{fig_optflow}c. The other high-order windows, Hann and Bump, exhibit a stationary flow as well. Intuitively, this hints that the optimization using higher order windows is more effective than the optimization using a Square-window, since the reduced vorticity of the stationary flow in Fig.~\ref{fig_optflow}c induces less turbulent viscosity, which finally reduces the effective drag of the airfoil.

A more detailed report of the different optimization runs is given in Fig.~\ref{fig_dsgn_res}. Here we can quantify the intuitive result of Fig.~\ref{fig_optflow}. The figure displays values of the optimization objective ${J}_w$ and the optimization constraint ${C}_w$ as well as the bound of the inequality constraint. Both lift ${C}_w$ and the lower bound of the lift constraint  ${C}_w=0.96$ are measured on the right hand side y-axis of Fig.~\ref{fig_dsgn_res}a. Note however, that the axis of Fig.~\ref{fig_dsgn_res}a is scaled by a factor $2$, since the values of ${J}_w$ and ${C}_w$ oscillate much more in the Square-window run compared to the other optimization runs. We can see, that the Square-windowed optimization struggles to satisfy the inequality constraint while minimizing the objective ${J}_w$. Each time the optimizer reduces ${J}_w$, it oversteps the lift constraint ${C}_w\geq0.94$, which results in an infeasible design. The subsequent design then increases both ${C}_w$ and ${J}_w$ but reduces ${J}_w$ a bit compared to the first design. This results in a very slow optimization procedure compared to the ones, where higher order windows are used.

Contrary to the Square-window run, the design process of a higher order window produces designs that do not repeat themselves in such a pattern, see Fig.~\ref{fig_dsgn_res}b-d. After two infeasible designs, the optimizer manages to reduce the windowed time averaged drag to  $30\%$ of the baseline value, whereas keeping the design feasible most of the times. 
These effects are explainable by considering the surface sensitivities in Subsection~\ref{subsec_surfsens}. We have seen, that the surface sensitivities computed with the Square-window change, when the time frame, over which the sensitivity was averaged, was shifted. 

During the shape optimization process, we change the airfoil design and hence, the limit cycle changes as well. This may result in a shifted phase, shape or period length of the new  {limit cycle}, compared to the baseline  {limit cycle}, or a combination of these effects, see \cite{Wilkins}. We have seen, that a higher order window mitigates the effect of a shifted  {limit cycle}, which is equivalent to shifting the averaging time frame, better than the Square-window. i.e. it is more robust with respect to a shifted  {limit cycle}. Another aspect of the increased oscillation in the Square-window design process, compared to the higher order design processes, is given by  the time-dependent sensitivities computed with different windows shown in Fig.~\ref{fig_sensOSC}b. As discussed in Subsection~\ref{sec_valWND} the design variable sensitivities computed using a Square-window oscillate depending on the length of the averaging time, and in this high Reynolds number test case, the amplitude even increases. Higher order windows however dampen their oscillatory behavior as time increases.

Whereas a changed shape of a period is not important for any windowing, including the Square-window, a change in the period length affects convergence speed of the windowed averaging, since a windowed time average converges depending on the number of completed periods in the given time frame $M$.

\subsubsection*{Test Case $Re=10^3$}
We consider briefly the test case at a lower Reynolds number. We consider again the optimization problem in Eq.~\eqref{eq_ASP_EX}, but at $Re=10^3$, $n_{tr}=500$ and $N=684$. The lift constraint $C_w$ is set to be greater equals $0.76$, which is the windowed time averaged lift of the baseline geometry. All other flow properties are set to the same values as in the test case above.

\begin{figure}
\begin{minipage}{.5\textwidth}
  \centering
  \includegraphics[width=.9\linewidth]{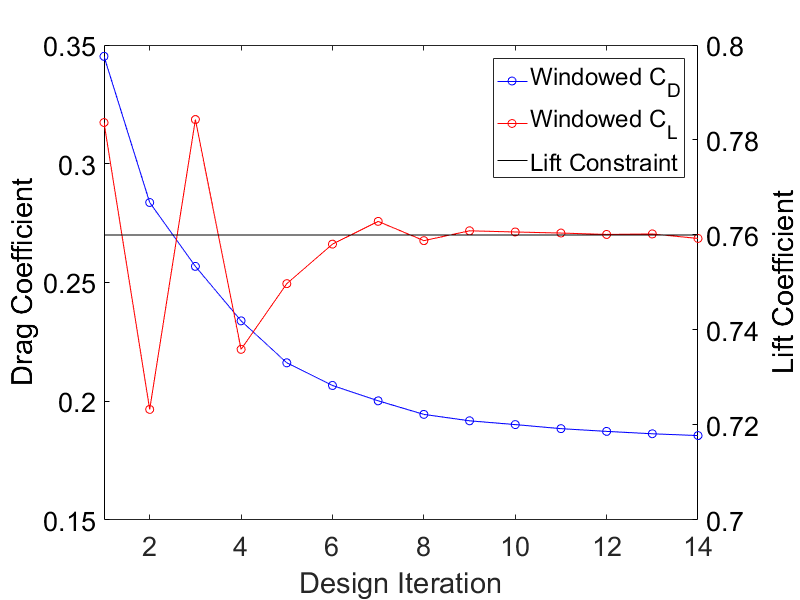}{\\(a) Square-windowing.}
\end{minipage}
\begin{minipage}{.5\textwidth}
  \centering
  \includegraphics[width=.9\linewidth]{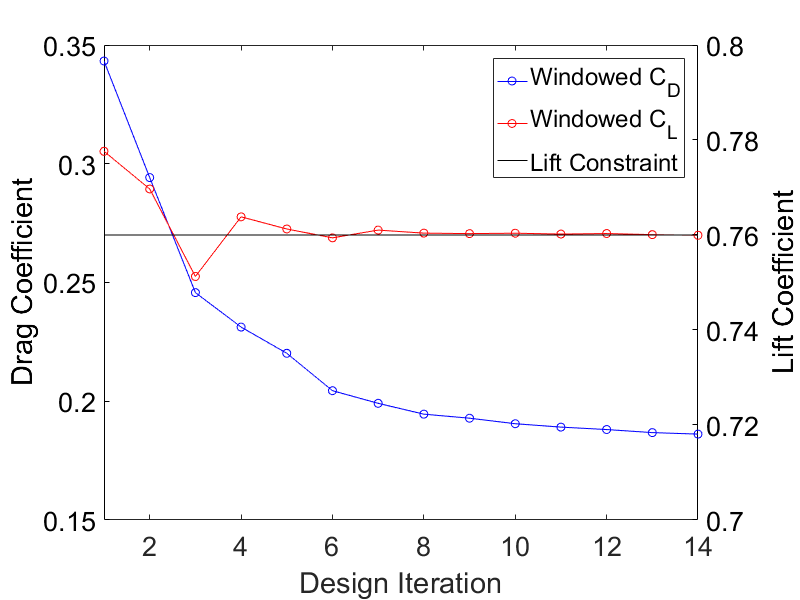}{\\(b) Hann-windowing.}
\end{minipage}
\begin{minipage}{.5\textwidth}
  \centering
  \includegraphics[width=.9\linewidth]{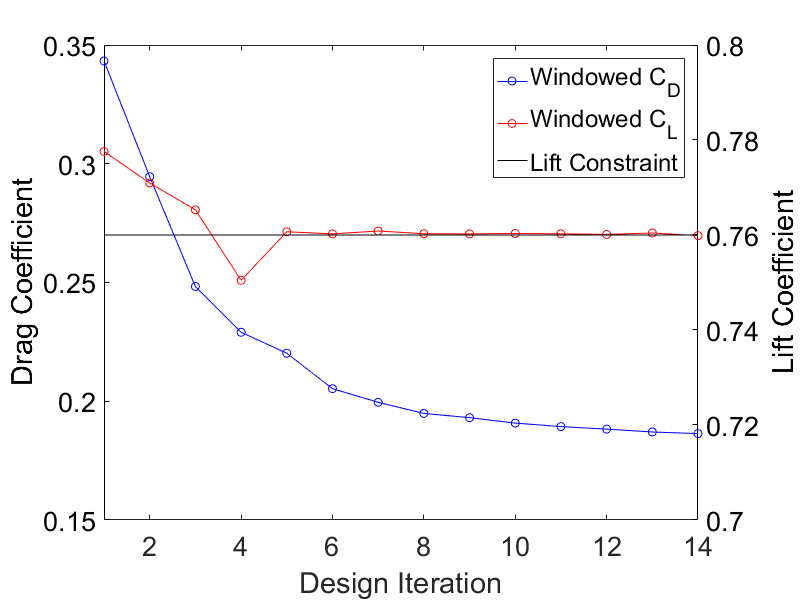}{\\(c) Hann-Square-windowing.}
\end{minipage}
\begin{minipage}{.5\textwidth}
  \centering
  \includegraphics[width=.9\linewidth]{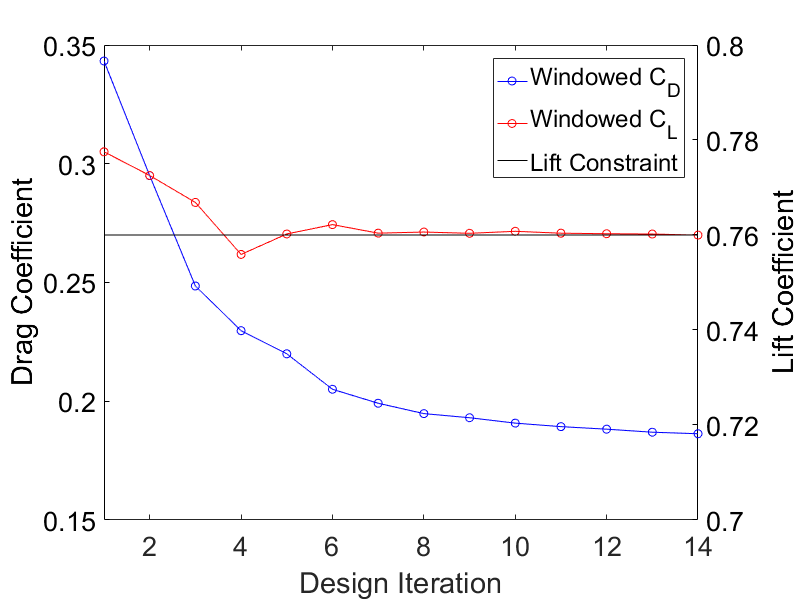}{\\(d)  Bump-windowing.}
\end{minipage}
 \caption{Shape optimization performed by SLSQP with different windows used. Both - feasible and infeasible - design steps are displayed. $Re=10^3$.}
 \label{fig_dsgn_resRe3}
\end{figure}
We have chosen a smaller time frame of $168$ iterations, which corresponds to approximately $5$ periods. We can see in Fig.~\ref{fig_sensOSCRe3}, that after $168$ iterations, the Square-windowed time average still oscillates, whereas higher order windowed time averages are already leveled.
Therefore, we expect worse results in case of a Square-windowed optimization compared to its higher order counterparts, which should perform comparably well.
Indeed, we can see in Fig.~\ref{fig_dsgn_resRe3}a, that the Square-windowed lift constraint $C_w$, which is colored red, oscillates for the first few design iterations and stays infeasible, i.e. below the lift constraint, for much more iterations than the higher order windowed lift constraints. The last infeasible Square-windowed lift appears at iteration $8$, whereas the last infeasible Hann-windowed value appears at iteration $6$, respectively $4$ in the case of Hann-Square- and Bump-windowing. 

The optimization objective, the windowed time averaged drag coefficient $J_w$, does not depend much on the chosen windowing function, as opposed to the results of the $Re=10^6$ test case. Furthermore, the oscillating behavior of the Square-windowed designs is not as extreme as the the Square-windowed designs of the $Re=10^6$ test case. 
This can be explained by considering the amplitude of both the instantaneous output as well as the windowed time averages of both test cases.
The amplitudes of the $Re=10^3$ test case at iteration $684$ is smaller by an order of magnitude than the corresponding results at $Re=10^6$ at $N=2200$.
Hence, the   {SLSQP} optimizer can handle the Square-windowed time averaged sensitivities of the $Re=10^3$ test case better than the ones of the $Re=10^6$ test case.
\newpage
\section{Periodically pitching Airfoil}
In this test case we consider NACA64A010 airfoil, modeled by a grid consisting of 12897 quadrilaterals and 23205 triangles. The airfoil surface consists of $250$ wall-boundary elements and at the far-field we have $68$ elements. We compute sensitivities with respect to $50$ design variables, which model the airfoil surface using Hicks-Henne functions.
The flow is computed with the URANS solver of SU2, where we use the Roe method for the Navier-Stokes fluxes and the Upwind method for the transport equation in the Spalart-Allmaras turbulence model. 
The angle of attack in this configuration is at $5$ degrees relatively low, and we achieve a periodic flow by tilting the wing around its longitudinal axis. Therefore, the airfoil increases or decreases its angle of attack. Here, we use a pitching amplitude of $3$ degrees and an angular frequency of $53.3491$ $\text{rad}/s$. The remaining free-stream parameters are given by
\begin{itemize}
 \item Mach Number: $0.35$
 \item Freestream Temperature: $88.15 K$
 \item Reynolds Number: $10^6$
 \item Reynolds Length: $1.0$ (Length of the airfoil)
\end{itemize} 
The simulation is computed with time step $\Delta t = 0.001636s$, which is equal to $72$ time steps per pitching period.

\subsection{Comparison of windowed time averages and windowed time averaged sensitivities}
\begin{figure}
\begin{minipage}{.5\textwidth}
  \centering
  \includegraphics[width=.9\linewidth]{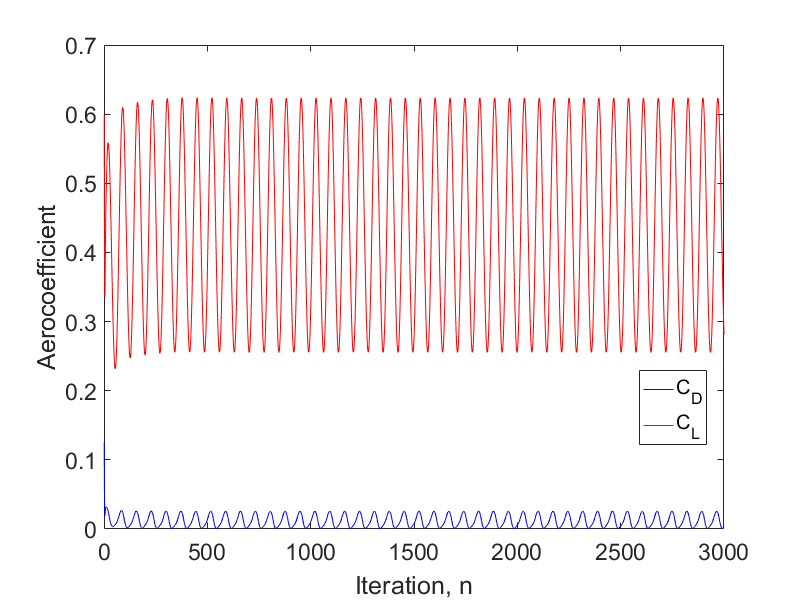}
\end{minipage}
\begin{minipage}{.5\textwidth}
  \centering
  \includegraphics[width=.9\linewidth]{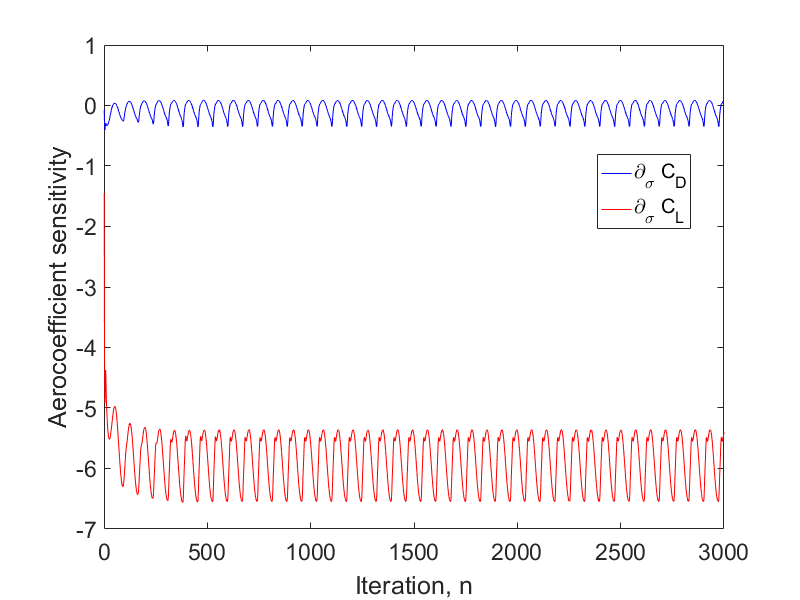}
\end{minipage}
\caption{Drag and Lift and their sensitivities of the pitching airfoil  over iteration number.}
\label{fig_pitchCDCL}
\end{figure}

Here, we display first the time-dependent drag and lift coefficient as well as their sensitivities in Fig.~\ref{fig_pitchCDCL}. We have a transient phase of approximately $130$ iterations and after it, both drag $C_D$ and lift $C_L$ oscillate with a period amplitude of approximately $0.19$ in the case of $C_D$, respectively $0.01$ in the case of $C_L$.
 We can see that, in contrast to Fig.~\ref{fig_RE}, which displays the NACA0012 test case at $Re=10^3$, the sensitivities of neither drag nor lift grow in amplitude. Their amplitudes are $0.1$, respectively $0.5$.
The effect of a non-increasing amplitude of the sensitivity can be explained by the fact, that the period length is only dependent on the angular frequency, which is independent of the design parameters. 

\begin{figure}
\begin{minipage}{.5\textwidth}
  \centering
  \includegraphics[width=.9\linewidth]{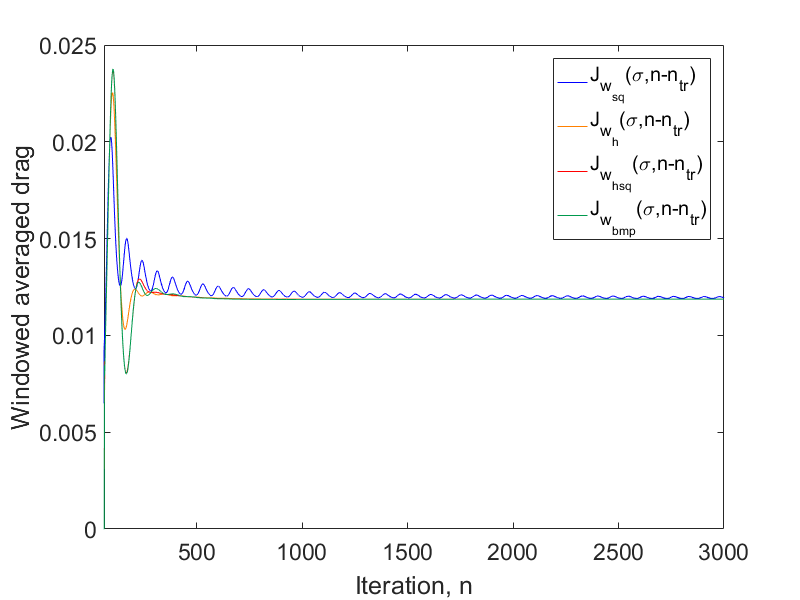}{\\(a) Windowed time average.}
\end{minipage}
\begin{minipage}{.5\textwidth}
  \centering
  \includegraphics[width=.9\linewidth]{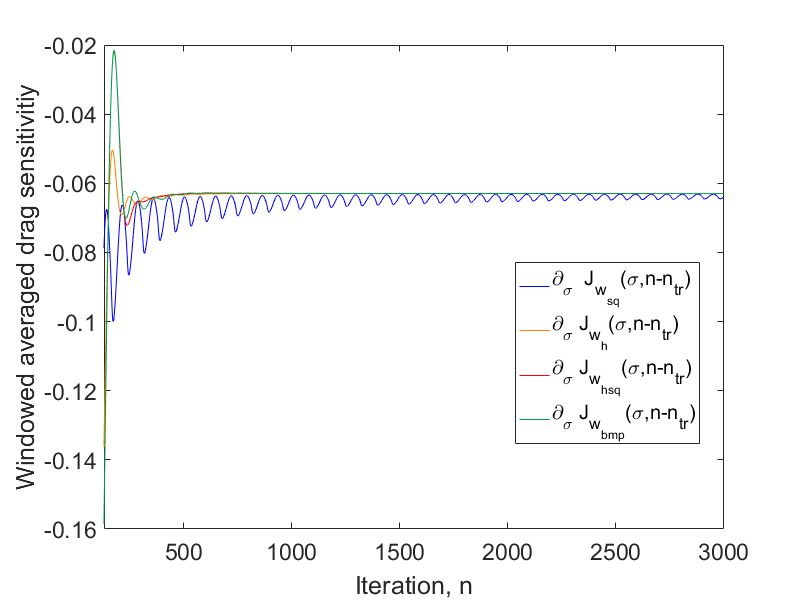}{\\(b)  Windowed time averaged sensitivity.}
\end{minipage}
 \caption{${J}_w(\sigma,n-n_{tr})$ and $\frac{\intD}{\intD \sigma} {J}_w(\sigma,n-n_{tr})$ with different windows over $n$, $n_{tr}=130$.}
\label{fig_wndPitch}
\end{figure}

We can see in Fig.~\ref{fig_wndPitch} the windowed time averages and their sensitivities. We start the windowed time average at $n_tr=130$, where the transient phase has passed. Both the windowed time averages and the windowed time averaged sensitivities converge for all higher order windows after a  couple of iterations. Furthermore, we can see that the Square-window performs almost equally well in case of the time averaged sensitivity as the time average. Reason for this is, that the period length is independent of the design parameter. This means that in Theorem~\ref{theo1}, we can use the convergence boundary given by Eq.~\eqref{eq_convNonSens} rather than Eq.~\eqref{eq_convSens} for the time averaged sensitivities, since the second equality in Eq.~\eqref{eq_avgSens1} holds in this case.

\begin{figure}
\begin{minipage}{.5\textwidth}
  \centering
  \includegraphics[width=.9\linewidth]{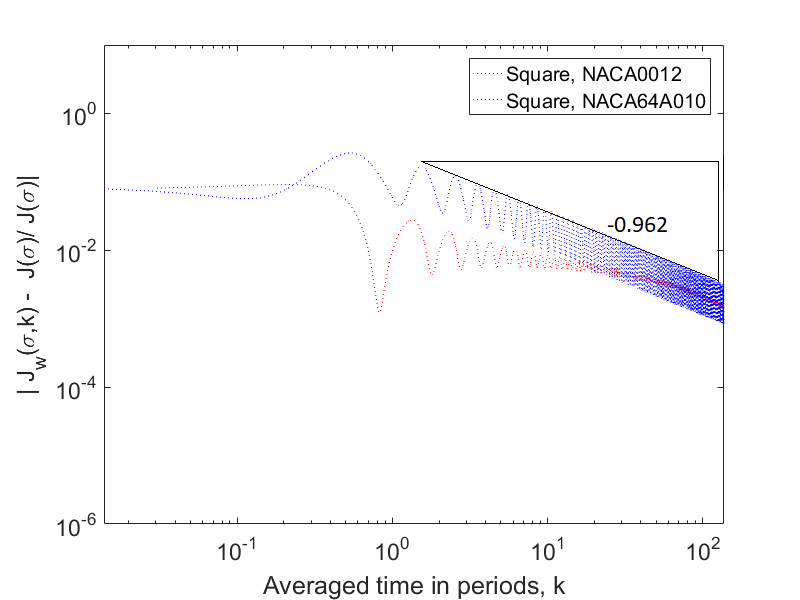}{\\(a) Windowed time average.}
\end{minipage}
\begin{minipage}{.5\textwidth}
  \centering
  \includegraphics[width=.9\linewidth]{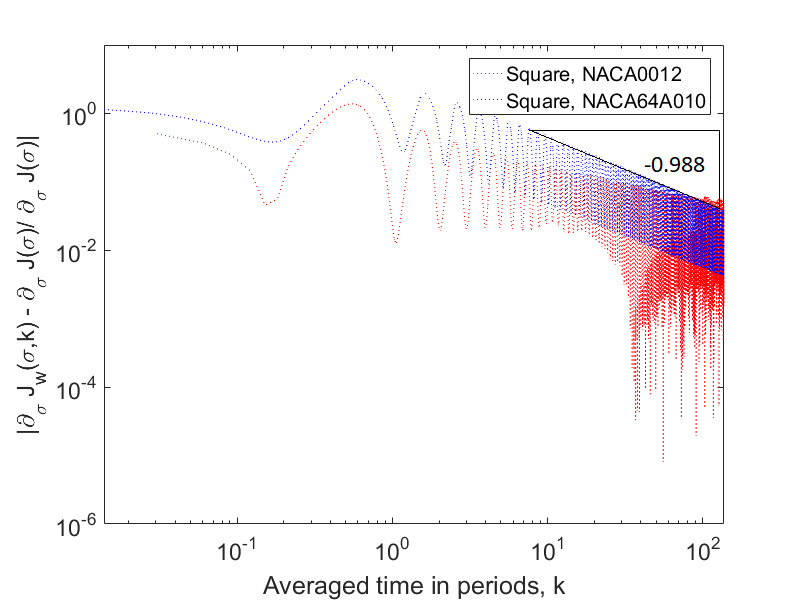}{\\(b) Windowed time averaged sensitivity.}
\end{minipage}
\caption{ Relative differences of Square-windows to their limits, where $k =\frac{N-n_{tr}}{T}$. We display the NACA0012 airfoil at $Re=10^3$ and the NACA64A10 airfoil.}
\label{fig_logCompare}
\end{figure}

Figure~\ref{fig_logCompare}a compares the Square-windowed time averages of the NACA0012 test case computed at $Re=10^3$ to the current NACA64A10 test case. Figure~\ref{fig_logCompare}b compares the corresponding Square-windowed time averaged sensitivities. The time scale of both figures is given by $k = \frac{N-n_{tr}}{T}$, where $T$ is the period length of limit cycle of the specific test case. The curves display the relative difference of the windowed time averages, respectively their sensitivities, to their corresponding limit values. Since a computation of the analytic value ${J}(\sigma)$ is not possible, use the Bump-windowed value for a very high value of $k$ as an approximation.    One can see that the Square-window in the NACA64A10 test case has approximately the same slope as its sensitivity, which confirms the consideration above. Note that the displayed value of the slope is corrected by the ratio of scales of the x- and y-axis. The slope of the windowed time average in the NACA0012 case asymptotically matches the counterpart of the NACA64A10 case. However, the Square-windowed time averaged sensitivity of the NACA64A10 case has a much steeper slope than its  NACA0012 counterpart.
This validates the consideration above. 

With this test case we have seen an example for which the method is not absolutely necessary. However, usage of higher order windows still speeds up the convergence rate of the time averaged output and the time averaged sensitivity.

\newpage
\section{Conclusion}
In this paper we have combined the unsteady aerodynamic optimization framework of SU2 with the windowing approach to obtain meaningful and robust design sensitivities. First we have defined the windowed time average as an output to the direct flow solver. Then we have embedded the long time windowing approach in the discrete adjoint solver of SU2 using a Lagrangian method. The discrete adjoint directly inherits the convergence properties of the primal flow solver.

Both the primal and adjoint solvers with windowed output functions have been applied in the NACA0012 and NACA64A010 test cases, which display an unsteady turbulent detached flow that exhibits limit cycle oscillations.
We have shown the superiority of the windowing-approach using high order windows, compared to the traditional non-windowed approach for sensitivity analysis of a period-average. Whereas non-windowed (or Square-windowed) averaged sensitivities tend to oscillate even for long time averages, high order windows quickly converge. We have demonstrated an increased robustness of computed sensitivities and a resulting more efficient optimization procedure using high order windows for the sensitivity calculation, compared to the traditional approach by taking the example of the NACA0012 airfoil as a two dimensional test case. Furthermore, high order windowing does not increase the computational cost of the sensitivity calculation.

\bibliography{sample}

\end{document}